\documentclass[12pt,centertags]{amsart}
\usepackage{euscript,amsmath,amstext,amsthm,a4,amssymb}
\usepackage{mathrsfs}
\usepackage{epsf}
\numberwithin{equation}{section}

\newcommand{\field}[1]{\mathbb{#1}}
\newcommand{\bZ}{\field{Z}}
\newcommand{\bR}{\field{R}}
\newcommand{\bC}{\field{C}}
\newcommand{\bN}{\field{N}}
\def\cC{\mathscr {C}}
\def\cL{\mathscr {L}}
\def\cO{\mathscr {O}}
\def\cR{\mathscr {R}}
\def\mA{\mathcal{A}}
\def\mB{\mathcal{B}}
\def\mC{\mathcal{C}}
\def\mL{\mathcal{L}}
\def\mO{\mathcal{O}}
\def\mQ{\mathcal{Q}}
\def\mR{\mathcal{R}}
\def\bE{{\bf E}}
\def\bk{{\bf k}}

\def\br{{\bf r}}

\def\Im{{\rm Im}}

\DeclareMathOperator{\End}{End}

\DeclareMathOperator{\Ker}{Ker}

\DeclareMathOperator{\Id}{Id}
\DeclareMathOperator{\supp}{supp}
\DeclareMathOperator{\tr}{Tr}

\newcommand{\spin}{$\text{spin}^c$ }
\newcommand{\spec}{\rm Spec}

\newcommand{\norm}[1]{\lVert#1\rVert}

\newcommand{\om}{\omega}

\newtheorem{thm}{Theorem}[section]
\newtheorem{lemma}[thm]{Lemma}
\newtheorem{prop}[thm]{Proposition}
\newtheorem{cor}[thm]{Corollary}
\theoremstyle{definition}
\newtheorem{defn}[thm]{Definition}
\theoremstyle{remark}
\newtheorem{rem}[thm]{Remark}
\newcommand{\be}{\begin{eqnarray}}
\newcommand{\ee}{\end{eqnarray}}
\newcommand{\wi}{\widetilde}
\newcommand{\var}{\varepsilon}

\newcommand{\comment}[1]{}

\begin{document}
\title{On the asymptotic expansion of Bergman kernel}
% \author{Xianzhe Dai, Kefeng Liu, Xiaonan Ma and Xiaowei Wang}

\author{Xianzhe Dai}
\address{Department of Mathematics, UCSB, CA 93106 USA (dai@math.ucsb.edu)}

\author{Kefeng Liu}
\address{Center of Mathematical Science, Zhejiang University and Department of Mathematics, UCLA, CA 90095-1555,
USA (liu@math.ucla.edu)}

%    Information for third author
\author{Xiaonan Ma}
\address{Centre de Math\' ematiques Laurent Schwartz, UMR 7640 du CNRS,
Ecole Polytechnique, 91128 Palaiseau Cedex,
France (ma@math.polytechnique.fr)}

%\thanks{\footnote{Department of Mathematics, UCLA,CA 90095-1555,
%USA (liu@math.ucla.edu)}\footnote{CNRS UMR 7640, Centre de
%Math\' ematiques, Ecole Polytechnique, 91128 Palaiseau Cedex,
%France (ma@math.polytechnique.fr)}}

\begin{abstract}
We study the asymptotic of the Bergman kernel of the \spin Dirac operator
on high tensor powers of a line bundle.
\end{abstract}

\maketitle

\section{Introduction} \label{s1}

The Bergman kernel in the context of several complex variables
(i.e. for pseudoconvex domains) has long been an important subject
(cf, for example, \cite{BFG}). Its analogue for complex projective
manifolds is studied in \cite{Tian}, \cite{Ru}, \cite{Zelditch},
\cite{Catlin}, \cite{Lu}, establishing the diagonal asymptotic expansion
for high powers of an ample line bundle. Moreover, the
coefficients in the asymptotic expansion encode geometric
information
 of the underlying complex projective manifolds. This asymptotic expansion
plays a crucial role in the recent work of \cite{D} where the
existence of K\"ahler metrics with constant scalar curvature is
shown to be closely related to Chow-Mumford stability.

Borthwick and Uribe \cite{BU1},
 Shiffman and Zelditch \cite{SZ02} were the first ones to study
the corresponding symplectic versions.
%the almost holomorphic sections of ample line bundles
% on symplectic manifolds.
Note that they use the almost holomorphic sections
 based on a construction of Boutet de
Monvel--Guillemin \cite{BoG} of a first order pseudodifferential operator
$D_b$ associated to the line bundle $L$ on a compact symplectic manifold,
 which mimic the ${\overline\partial}_b$ operator on the circle bundle
in the holomorphic case.
 The Szeg{\"o} kernels are well defined
modulo  smooth operators on the associated circle bundle,
 even though $D_b$ is neither canonically defined nor unique
(Actually, Boutet de Monvel--Guillemin define first the Szeg\"o kernels, 
then construct the operator $D_p$ from the Szeg\"o kernels).
 Moreover, in the holomorphic case, the Szeg{\"o} kernels are exactly
(modulo  smooth operators)
the Szeg\"o kernel associated to the holomorphic sections
by Boutet de Monvel-Sj\"ostrand \cite{BS}.
In the very important paper \cite{SZ02}, Shiffman and Zelditch also gave
a simple way to construct first the Szeg{\"o} kernels, then the
operator $D_b$ from the construction of
 Boutet de Monvel--Guillemin \cite{BoG},
and in \cite[Theorem 1]{SZ02}, they studied the near diagonal
asymptotic expansion and small ball Gaussian estimate (for
$d(x,y)\leq C/\sqrt{p}$ where $p$ is the power of the line bundle $L$).
% In particular, they obtained the near diagonal asymptotic
%expansion of the Szeg{\"o} kernel in the holomorphic case.
 On the other hand, in the holomorphic setting,
in \cite{Christ91}, Christ (and  Lindholm in \cite{Li}) proved an
Agmon type estimate for the Szeg{\"o} kernel on $\bC^1$, but they
did not treat the asymptotic expansions.

In this paper, we establish the full off-diagonal asymptotic
expansion and Agmon estimate for the Bergman kernel of the
spin$^c$ Dirac operator associated to high powers of an ample line
bundle in the general context of symplectic manifolds and
orbifolds (Cf. Theorem \ref{tue17}; note the important factors on the right hand side of the
estimate \eqref{ue66} which make our estimate uniform for $Z$ and $Z'$). Our motivations are to extend Donaldson's
work \cite{D} to orbifolds and to understand the relationship
between heat kernel, index formula and stability. Moreover, the spin$^c$
Dirac operator is a natural geometric operator associated to the
symplectic structure. As a result the coefficients in the
asymptotic expansion are naturally polynomials of the curvatures
and their derivatives.

Let $(X,\om)$ be a compact symplectic manifold of real dimension
$2n$. Assume that there exists a Hermitian line bundle $L$ over
$X$ endowed with a Hermitian connection $\nabla^L$ with the
property that $\frac{\sqrt{-1}}{2\pi}R^L=\omega$, where
$R^L=(\nabla^L)^2$ is the curvature of $(L,\nabla^L)$. Let
$(E,h^E)$ be a Hermitian vector bundle  on $X$ with Hermitian
connection $\nabla^E$ and its curvature $R^E$.

Let $g^{TX}$ be a Riemannian metric on $X$.
Let ${\bf J}:TX\longrightarrow TX$ be the skew--adjoint linear
map which satisfies the relation
\be \label{0.1}
\om(u,v)=g^{TX}({\bf J}u,v)
\ee
 for $u,v \in TX$.
Let $J$ be an almost complex structure which is separately
compatible with $g^{TX}$ and $\om$, and $\om(\cdot,
J\cdot)$ defines a metric on $TX$.
 Then $J$ commutes with ${\bf J}$ and $-J{\bf J}\in \End(TX)$ is positive,
thus $-J{\bf J}=(-{\bf J}^2)^{1/2}$. 
Let $\nabla ^{TX}$ be the Levi-Civita connection on $(TX,
g^{TX})$ with curvature $R^{TX}$, and  $\nabla^{TX}$ induces a
natural connection $\nabla ^{\det}$ on $\det (T^{(1,0)} X)$ with
curvature $R^{\det}$ (cf. Section \ref{s2}). The \spin Dirac
operator $D_p$ acts on
$\Omega^{0,{\scriptscriptstyle{\bullet}}}(X,L^p\otimes E)
=\bigoplus_{q=0}^n\Omega^{0,q}(X,L^p\otimes E)$, the direct sum of
spaces of $(0,q)$--forms with values in $L^p\otimes E$.

Let $\{S^p_i\}_{i=1}^{d_p}$
$(d_p = \dim \Ker D_p)$ be any orthonormal basis of
$\Ker D_p$ with respect to the inner
 product (\ref{b3}). We define
the diagonal of the  Bergman kernel of $D_p$  (the distortion
function) by
\begin{align} \label{0.2}
&B_p(x) = \sum_{i=1}^{d_p} S^p_i (x) \otimes (S^p_i(x))^*
\in \End (\Lambda (T^{*(0,1)}X)\otimes E)_x
\end{align}
Clearly $B_p(x)$ does not depend on the choice of $\{S^p_i\}$. We
denote by $ I_{\bC\otimes E}$  the projection from $\Lambda
(T^{*(0,1)}X)\otimes E$ onto $\bC\otimes E$ under the decomposition
$\Lambda (T^{*(0,1)}X)= \bC \oplus \Lambda ^{>0} (T^{*(0,1)}X)$. Let
$\det {\bf J}$ be the determinant function of ${\bf J}_x\in
\End(T_xX)$, and  $|{\bf J}|= (-{\bf J}^2)^{1/2}\in \End(T_x X)$. A
simple corollary of Theorem \ref{tue17}$^\prime$ is:

\begin{thm}\label{t0.1} There exist smooth coefficients
$b_r(x)\in \End (\Lambda (T^{*(0,1)}X)\otimes E)_x$ which
are polynomials in $R^{TX}$, $R^{\det}$,
$R^E$ {\rm (}and $R^L${\rm )}  and their derivatives with order
$\leq 2r-1$  {\rm (}resp. $2r${\rm )} and reciprocals
of linear combinations of eigenvalues of ${\bf J}$  at $x$, and
 $b_0=(\det {\bf J})^{1/2} I_{\bC\otimes E}$,
such that for any $k,l\in \bN$, there exists
$C_{k,l}>0$ such that for any $x\in X$, $p\in \bN$,
\begin{align}\label{0.3}
&\Big |B_p(x)
- \sum_{r=0}^{k} b_r(x) p^{n-r} \Big |_{\cC^l} \leq C_{k,l} p^{n-k-1}.
\end{align}
Moreover,
the expansion is uniform
in that for any  $k,l\in \bN$, there is an integer $s$ such that
if all data {\rm(}$g^{TX}$, $h^L$, $\nabla ^L$, $h^E$, $\nabla ^E${\rm)}
run over a set which are bounded in $\cC^s$ and with $g^{TX}$
 bounded below, there exists the constant  $C_{k,\,l}$
independent of $g^{TX}$,
and the $\cC^l$-norm in \eqref{0.3} includes also the
derivatives on the parameters.
\end{thm}

We also study the asymptotic expansion of the corresponding heat
kernel and relates it to that of the Bergman kernel. Let
$\exp(-\frac{u}{p} D_p^2) (x,x')$ be the smooth kernel of
$\exp(-\frac{u}{p} D_p^2)$ with respect to the Riemannian volume
form $dv_X(x')$. We introduce in \eqref{b4}, $\om_d (x)\in
\End(\Lambda (T_x^{*(0,1)}X))$.

\begin{thm}\label{t0.2} There exist smooth sections
 $b_{r,u}$ of $\End (\Lambda (T^{*(0,1)}X)\otimes E)$ on $X$
which are polynomials in $R^{TX}$, $R^{\det}$,
$R^E$ {\rm (}and $R^L${\rm )}  and their derivatives with order
$\leq 2r-1$  {\rm (}resp. $2r${\rm )} and
 functions on the eigenvalues of ${\bf J}$  at $x$, and
 $b_{0,u}=\Big(\det (\frac{|{\bf J}|}{1- e ^{-4\pi u|{\bf J}|}}) \Big)^{1/2}
e ^{-4\pi u \om_d}$,
such that for each $u>0$ fixed, we have the asymptotic expansion
in the sense of \eqref{0.3} as $p\to \infty$,
\be\label{0.4}
\exp(-\frac{u}{p} D_p^2) (x,x) = \sum_{r=0}^k b_{r,u}(x) p^{n-r} +
\cO(p^{n-k-1}). \ee
Moreover, there exists $c>0$ such that as
$u\to +\infty$,
\be\label{0.5}
 b_{r,u}(x) = b_r(x) +\cO(e ^{-cu}).
\ee
\end{thm}
Note that the coefficient $b_{0,u}$ in Theorem \ref{t0.2} was first
obtained in \cite[(f)]{B}. Theorems \ref{t0.1}, \ref{t0.2} give
 us a way to compute the coefficient  $b_r(x)$,
as it is relatively easy
to compute $ b_{r,u}(x)$ (cf. (\ref{ue49}), (\ref{ue70})).
As an example, we compute $b_1$ which plays an important role in
Donaldson's recent work \cite{D}. Note
if  $(X,\omega)$ is K\"ahler and ${\bf J}=J$, then
$B_p(x)\in \cC^\infty(X, \End(E))$ for $p$ large enough,
thus $b_r(x)\in \End(E)_x$.

\begin{thm}\label{t0.3} If $(X,\omega)$ is K\"ahler and ${\bf J}=J$,
then there exist smooth functions $b_r(x)\in \End (E)_x$ such that
we have  \eqref{0.3}, and $b_r$ are polynomials in $R^{TX}$, $R^E$
and their derivatives with order $\leq 2r-1$  at $x$. Moreover, \be
\label{0.6} b_0=\Id _E, \quad b_1 = \frac{1}{4\pi} \Big [ \sqrt{-1}
\sum_i R^E (e_i, Je_i) +\frac{1}{2}r^X \Id_E \Big ]. \ee here $r^X$
is the scalar curvature of $(X, g^{TX})$, and $\{e_i\}$ is an
orthonormal basis of  $(X, g^{TX})$.
\end{thm}

Theorem \ref{t0.3} was essentially obtained in \cite{Lu},
\cite{Wang1} by applying the peak section trick, and in
\cite{Catlin}, \cite{Zelditch} and \cite{Charles} by applying the
Boutet de Monvel-Sj\"ostrand parametrix for the Szeg\"o kernel
\cite{BS}. We refer the reader to \cite{D}, \cite{Wang1} for its
interesting applications.

Our proof of Theorems \ref{t0.1}, \ref{t0.2} is inspired by local
Index Theory, especially by \cite[\S 11]{BL}, and we derive Theorem
\ref{t0.1} from Theorem \ref{t0.2}. In particular, with the help of
the heat kernel, we get the full off-diagonal asymptotic expansion
for the Bergman kernel and the Agmon estimate for the remainder term
of the asymptotic expansion (Cf. Theorem \ref{tue17}). And when
$(X,\omega)$ is a K\"ahler manifold, ${\bf J}=J$ on $X$ and $E=\bC$,
we recover \cite[Theorem 1]{SZ02} if we restrict
 Theorem \ref{tue17}$^\prime$ to  $|Z|, |Z'|<C/\sqrt{p}$.

One of the advantages of our method is that it can be easily
generalized to the orbifold situation, and indeed, in \eqref{0e9},
we deduce the explicit asymptotic expansion near the singular set of
the orbifold.

\begin{thm}  \label{t0.4} If $(X,\omega)$ is a symplectic orbifold
with the singular set $X'$, and $L$, $E$ are
corresponding proper orbifold vector bundles on $X$ as in Theorem \ref{t0.1}.
 Then there exist smooth coefficients
$b_r(x)\in \End (\Lambda (T^{*(0,1)}X)\otimes E)_x$
with $b_0=(\det {\bf J})^{1/2} I_{\bC\otimes E}$, and $b_r(x)$ are
polynomials in $R^{TX}$, $R^{\det}$,
$R^E$ {\rm (}and $R^L${\rm )}  and their derivatives with order
$\leq 2r-1$ {\rm (}resp. $2r${\rm )} and reciprocals
of linear combinations of eigenvalues of ${\bf J}$  at $x$,
such that for any $k,l\in \bN$, there exist
$C_{k,l}>0$, $N\in \bN$ such that for any $x\in X$, $p\in \bN$,
\begin{multline}\label{0.7}
\Big |\frac{1}{p^n}B_p(x)
- \sum_{r=0}^{k} b_r(x) p^{-r} \Big |_{\cC^l}
\leq C_{k,l}
\Big (p^{-k-1} + p^{l/2}(1+\sqrt{p}d(x,X') )^N
e^{-C \sqrt{p} d(x, X')}\Big ).
\end{multline}
Moreover if the orbifold $(X,\omega)$ is K\"ahler,  ${\bf J}=J$
and the proper orbifold vector bundles $E,L$ are holomorphic on $X$,
then  $b_r(x)\in \End (E)_x$ and  $b_r(x)$ are polynomials in $R^{TX}$,
$R^E$  and their derivatives with order
$\leq 2r-1$  at $x$.
\end{thm}

This paper is organized as follows. In Section \ref{s2}, we recall a
 result on the spectral gap of the \spin Dirac operator \cite{MM}.
  In  Section \ref{s3}, we localizes the problem by finite propagation speed
and use the rescaling in local index theorem to prove  Theorems \ref{t0.1}, \ref{t0.2}.
 In  Section \ref{s4}, we compute the coefficients of the asymptotic expansion and explain
how to generalize our method to the orbifold situation.

The results of this paper have been announced in \cite{DKM}.

%%%%%%%%%%%%%%%%%%%%%%%%%%%%%%%%%%%%%%%%%%%%%%%%%%%%%%%%%%%%%%%%%%%%%%%%%%%%%%%%%%%%%%%%%%%%%%%%%%%%%

\subsection*{Aknowledgements}
We thank Professors Jean-Michel Bismut,  Jean-Michel Bony and
Johannes Sj\" ostrand for useful conversations. We also thank
Xiaowei Wang for useful discussions, and Laurent Charles for
sending us his paper. Finally, the authors acknowledge useful
comments and suggestions from the referee which help improve the
paper and clarify the relationship with some previous work.
\section{The spectral gap of the \spin Dirac operator}\label{s2}

The almost complex structure $J$ induces a splitting
$T_\bR X\otimes_\bR \bC=T^{(1,0)}X\oplus T^{(0,1)}X$,
where $T^{(1,0)}X$ and $T^{(0,1)}X$
are the eigenbundles of $J$ corresponding to the eigenvalues $\sqrt{-1}$
and $-\sqrt{-1}$ respectively. Let $T^{*(1,0)}X$  and $T^{*(0,1)}X$
be the corresponding dual bundles.
For any $v\in TX$ with decomposition $v=v_{1,0}+v_{0,1}
\in T^{(1,0)}X\oplus T^{(0,1)}X$,  let ${\overline v^\ast_{1,0}}\in T^{*(0,1)}X$ be the metric
dual of $v_{1,0}$. Then $c(v)=\sqrt{2}({\overline v^\ast_{1,0}}\wedge-
i_{v_{\,0,1}})$ defines the Clifford action of $v$ on
$\Lambda (T^{*(0,1)}X)$, where $\wedge$ and $i$
denote the exterior and interior product respectively. Set
\begin{equation}\label{b1}
\mu_0=\displaystyle\inf_{{u\in T_x^{(1,0)}X,\,x\in X}} R^L_x(u,\overline{u})/|u|^2_{g^{TX}} >0.
\end{equation}

Let $\nabla^{TX}$ be the Levi-Civita connection
of the metric $g^{TX}$.
By \cite[pp.397--398]{LM}, $\nabla^{TX}$ induces canonically
a Clifford connection $\nabla^{\text{Cliff}}$ on $\Lambda (T^{*(0,1)}X)$
(cf. also \cite[\S 2]{MM}).
Let $\nabla^{E_p}$ be the connection on
$E_p=\Lambda (T^{*(0,1)}X)\otimes L^p\otimes E$
induced by $\nabla^{\text{Cliff}}$, $\nabla^L$ and  $\nabla^E$.

Let $\langle \quad \rangle_{E_p}$ be the metric on $E_p$ induced by
 $g^{TX}$, $h^L$ and $h^E$.
Let $dv_X$ be the Riemannian volume form of $(TX, g^{TX})$.
The $L^2$--scalar product on $\Omega^{0,\scriptscriptstyle{\bullet}}(X,L^p\otimes E)$, the space of smooth sections of
$E_p$, is given by
\begin{equation}\label{b3}
\langle s_1,s_2 \rangle =\int_X\langle s_1(x),
s_2(x)\rangle_{E_p}\,dv_X(x)\,.
\end{equation}
We denote the corresponding norm with $\norm{\cdot}_{L^ 2}$.
Let $\{e_i\}_i$ be an orthonormal basis of $TX$.

\begin{defn}\label{Dirac}
The \spin Dirac operator $D_p$ is defined by
\begin{equation}\label{defDirac}
D_p=\sum_{j=1}^{2n}c(e_j)\nabla^{E_p}_{e_j}:
\Omega^{0,\scriptscriptstyle{\bullet}}(X,L^p\otimes E)\longrightarrow
\Omega^{0,\scriptscriptstyle{\bullet}}(X,L^p\otimes E)\,.
\end{equation}
\end{defn}
\noindent
$D_p$ is a formally self--adjoint, first order elliptic differential operator on
$\Omega^{0,\scriptscriptstyle{\bullet}}(X,L^p\otimes{E})$,
which interchanges $\Omega^{0,\text{even}}(X,L^p\otimes E)$
and $\Omega^{0,\text{odd}}(X,L^p\otimes E)$.

We denote by  $P^{T^{(1,0)}X}$
  the projection from $T_\bR X\otimes_\bR \bC$ to $T^{(1,0)}X$.
Let $\nabla^{T^{(1,0)}X}=P^{T^{(1,0)}X}\,\nabla^{TX}P^{T^{(1,0)}X}$
 be the Hermitian connection on $T^{(1,0)}X$ induced by  $\nabla^{TX}$
with curvature  $R^{T^{(1,0)}X}$. Let $\nabla^{\det (T^{(1,0)}X)}$
be the connection on $\det (T^{(1,0)}X)$ induced by  $\nabla^{T^{(1,0)}X}$
with curvature $R^{\det}=\tr[ R^{T^{(1,0)}X}]$.
Let $\{w_i\}$ be an orthonormal frame of $(T^{(1,0)}X, g^{TX})$.
Set
\begin{align}  \label{b4}
&\om_d=-\sum_{l,m} R^L (w_l,\overline{w}_m)\,\overline{w}^m\wedge
\,i_{\overline{w}_l}\,,
\qquad \tau(x)=\sum_j R^L (w_j,\overline{w}_j)\,.
\end{align}
Let $r^X$ be the  scalar curvature of $(TX,g^{TX})$, and
\begin{equation*}
\mathbf{c}(R)=\sum_{l<m}\left(R^E+\tfrac{1}{2}\tr\left[R^{T^{(1,0)}X}\right]
\right)(e_l,e_m)\,c(e_l)\,c(e_m)\,.
\end{equation*}
Then the Lichnerowicz formula \cite[Theorem 3.52]{BeGeV}
(cf. \cite[Theorem 2.2]{MM}) for $D_p^2$ is
\begin{equation}\label{Lich}
D^2_p=\left(\nabla^{E_p}\right)^{\ast}\,
\nabla^{E_p}-2p\om_d-
p\tau+\tfrac{1}{4}r^X+ \mathbf{c}(R),
\end{equation}

If $A$ is any operator, we denote by $\spec (A)$ the spectrum of $A$.

The following simple result was obtained in \cite[Theorems 1.1, 2.5]{MM}
by applying the Lichnerowicz formula
(cf. also \cite[Theorem 1]{BiV} in the holomorphic case).

\begin{thm}\label{t2.1}
There exists $C_L>0$ such that for any $p\in \bN$ and
any $s\in\Omega^{>0}(X,L^p\otimes E)
=\bigoplus_{q\geqslant 1}\Omega^{0,q}(X,L^p\otimes E)$,
\begin{equation}\label{main1}
\norm{D_{p}s}^2_{L^ 2}\geqslant(2p\mu_0-C_L)\norm{s}^2_{L^ 2}\,.
\end{equation}
Moreover  $\spec D^2_p \subset \{0\}\cup [2p\mu_0 -C_L,+\infty[$.
\end{thm}

\section{Bergman kernel}\label{s3}

In this Section, we will study the uniform estimate with
its derivatives on $t=\frac{1}{\sqrt{p}}$ of the heat kernel
and the Bergman kernel of $D^2_p$ as $p\to \infty$.
The first difficulty is that the space $\Omega^{0,\bullet}(X,L^p\otimes E)$
 depends on $p$.
To overcome this, we will localize the problem to a problem on
 $\bR^{2n}$. Now, after rescaling, another substantial difficulty appears, which
is the lack of the usual elliptic estimate on $\bR ^{2n}$ for the
rescaled Dirac operator. Thus we introduce a family of Sobolev
norms defined by the rescaled connection on $L^p$, then we can
extend  the functional analysis technique developed in \cite[\S
11]{BL}, and in this way, we can even get the estimate on its
derivatives on  $t=\frac{1}{\sqrt{p}}$.

This section is organized as follows. In Section \ref{s3.1}, we
establish the fact that the asymptotic expansion of $B_p(x)$ is
local on $X$. In Section \ref{s3.2}, we derive an asymptotic
expansion of $D_p$ in normal coordinate. In Section \ref{s3.3}, we
study the uniform estimate with its derivatives on $t$ of the heat
kernel and the Bergman kernel  associated to the rescaled operator
$L^t_2$ from $D_p^2$. In Theorem \ref{tue14}, we estimate
uniformly the remainder term of the Taylor expansion of  $e
^{-uL^t_2}$
 for $u\geq u_0>0,t\in [0,1]$. In Section \ref{s3.5}, we identify $J_{r,u}$
the coefficient of the Taylor expansion of  $e ^{-uL^t_2}$ with
the Volterra expansion of the heat kernel, thus giving us a way to
compute the coefficient $b_j$ in  Theorem \ref{t0.1}. In Section
\ref{s3.5},  we prove Theorems \ref{t0.1}, \ref{t0.2}.

\subsection{Localization of the problem}\label{s3.1}

Let $a^X$ be the injectivity radius of $(X, g^{TX})$,
and $\var\in (0,a^X/4)$.
We denote by $B^{X}(x,\epsilon)$ and  $B^{T_xX}(0,\epsilon)$ the open balls
in $X$ and $T_x X$ with center $x$ and radius $\epsilon$, respectively.
Then the map $ T_x X\ni Z \to \exp^X_x(Z)\in X$ is a
diffeomorphism from $B^{T_xX}(0,\epsilon)$  on $B^{X}(x,\epsilon)$ for
$\epsilon \leq a^X$.  From now on, we identify $B^{T_xX}(0,\epsilon)$
with $B^{X}(x,\epsilon)$ for $\epsilon\leq a ^X$.

Let $f : \bR \to [0,1]$ be a smooth even function such that
 \begin{align} \label{c2}
f(v) = \left \{ \begin{array}{ll}  1 \quad {\rm for}
\quad |v| \leq  \var/2, \\
  0 \quad {\rm for} \quad |v| \geq  \var.
\end{array}\right.
\end{align}
Set
 \begin{align} \label{1c2}
F(a)= \Big(\int_{-\infty}^{+\infty}f(v) dv\Big)^{-1} \int_{-\infty}^{+\infty} e ^{i v a} f(v) dv.
\end{align}
Then $F(a)$ lies in Schwartz space $\mathcal{S} (\bR)$ and $F(0)=1$.

Let $P_p$ be the orthogonal projection from
$\Omega^{0,\bullet}(X,L^p\otimes E)$
on $\Ker D_p$, and let  $P_p(x,x')$, $F(D_p)(x,x')$ ($x,x'\in X$),
 be the smooth kernels of $P_p$, $F(D_p)$
with respect to the volume form $dv_X(x')$.
The kernel $P_p(x,x')$ is called the Bergman kernel of $D_p$.
 By (\ref{0.2}),
\be\label{c1}
B_p(x)=P_p(x,x).
\ee

\begin{prop}\label{0t3.0}
 For any $l,m\in \bN$, $\var>0$, there exists $C_{l,m,\var}>0$
such that for $p\geq 1$, $x,x'\in X$,
\begin{align}\label{1c3}
&|F(D_p)(x,x') - P_p(x,x')|_{\cC^m(X\times X)} \leq C_{l,m,\var} \, p^{-l},\\
&|P_p(x,x')|_{\cC^m(X\times X)} \leq C_{l,m,\var}\, p^{-l} \quad \mathrm{if} \,
d(x,x') \geq \var \nonumber.
\end{align}
Here the $\cC^m$ norm is induced by $\nabla^L, \nabla^E$
and $\nabla^{\text{Cliff}}$.
\end{prop}

\begin{proof} For $a\in \bR$, set
\begin{eqnarray}\label{0c2}
\phi_p(a) = 1_{[\sqrt{p\mu_0}, +\infty[} (|a|) F(a).
\end{eqnarray}
Then by Theorem \ref{t2.1}, for $p> C_L/\mu_0$,
\begin{eqnarray}\label{0c3}
F(D_p)-P_p = \phi_p(D_p).
\end{eqnarray}
By (\ref{1c2}), for any $m\in \bN$ there exists $C_{m}>0$
such that
\be\label{1c9}
\sup_{a\in \bR} |a|^m |F(a) | \leq C_{m}.
\ee

As $X$ is compact, there exist $\{x_i\}_{i=1}^r$ such that
 $\{U_i = B^X(x_i,\var)\}_{i=1}^r$ is a covering of $X$.
We identify $B^{T_{x_i}X}(0,\var)$ with $B^{X} (x_i,\var)$
by geodesic as above. We identify
$(TX)_Z, (E_p)_Z$ for $Z\in B^{T_{x_i}X}(0,\var)$
to $T_{x_i}X, (E_p)_{x_i}$  by parallel transport with respect
to the connections $\nabla ^{TX}$, $\nabla^{E_p}$ along the curve
$\gamma_Z: [0,1]\ni u \to \exp^X_{x_i} (uZ)$. Let $\{e_i\}_i$
be an orthonormal basis of $T_{x_i}X$. Let $\wi{e}_i (Z)$ be the parallel
transport of ${e}_i$ with respect to $\nabla^{TX}$ along the above curve.
Let $\Gamma ^E, \Gamma ^L, \Gamma ^{\text{Cliff}}$
be the corresponding connection forms of $\nabla^E$, $\nabla^L $
and $\nabla^{\text{Cliff}}$ with respect to any fixed frame for
$E,L$, $\Lambda (T^{*(0,1)}X)$ which is parallel along the curve $\gamma_Z$
under the trivialization on $U_i$.

Denote by  $\nabla_U$  the ordinary differentiation
 operator on $T_{x_i}X$ in the direction $U$.
Then
\be\label{c10}
D_p =\sum_{j} c(\wi{e}_j) \Big ( \nabla_{\wi{e}_j} +p \Gamma ^L(\wi{e}_j)
+ \Gamma ^{\text{Cliff}}(\wi{e}_j)+ \Gamma ^E(\wi{e}_j)  \Big ).
\ee
Let $\{ \varphi_i \}$ be a partition of unity subordinate to $\{U_i\}$.
For $l\in \bN$, we define a Sobolev norm on the $l$-th Sobolev space
$H^l(X,E_p)$ by
\be\label{c11}
\| s\| _{H^l_p}^2 = \sum_i \sum_{k=0}^l \sum_{i_1 \cdots i_k=1} ^{2n}
\|\nabla_{e_{i_1}}\cdots  \nabla_{e_{i_k}}(\varphi _i s)\|_{L^2}^2
\ee
Then by (\ref{c10}), there exists $C>0$ such that for $p\geq 1$, $s\in
H^1(X, E_p)$,
\be\label{c12}
\|s\|_{H^1_p}   \leq C(\|D_p s\|_{L^2} + p\|s\|_{L^2}).
\ee
Let $Q$ be a differential operator of order $m\in \bN$ with
scalar principal symbol and with  compact support
in $U_i$, then
\be \label{c14} \qquad
[D_p,Q] = \sum_{j}  p [c(\wi{e}_j) \Gamma ^L(\wi{e}_j), Q]
+\sum_{j}  \Big  [c(\wi{e}_j)\Big  (\nabla_{\wi{e}_j}
+ \Gamma ^{\text{Cliff}}(\wi{e}_j)+ \Gamma ^E(\wi{e}_j)\Big ),  Q\Big ]
\ee
which are differential operators  of order $m-1$, m respectively.
By (\ref{c12}), (\ref{c14}),
\begin{align}\label{c15}
\|Qs\|_{H^1_p} &  \leq C(\|D_pQ s\|_{L^2} + p\|Qs\|_{L^2})\\
 & \leq C(\|Q D_p s\|_{L^2} + p\|s\|_{H^m_p}).\nonumber
\end{align}
From (\ref{c15}), for $m\in \bN$, there exists $C'_m>0$ such that for
$p\geq 1$,
\be\label{c16}
\|s\|_{H^{m+1}_p} \leq C'_m  (\|D_p s\|_{H^m_p} + p \|s\|_{H^m_p}).
\ee
This means
\be\label{c17}
\|s\|_{H^{m+1}_p} \leq C'_m \sum_{j=0}^{m+1}  p^{m+1-j}\|D_p^j s\|_{L^2}.
\ee
Moreover from
$\langle D_p^{m'} \phi_p(D_p)Q s,s'\rangle
=\langle s,Q^ * \phi_p(D_p) D_p^{m'} s'\rangle$,
(\ref{0c2}) and (\ref{1c9}), we know that for
$l,m'\in \bN$, there exists $C_{l,m'}>0$ such that for $p\geq 1$,
\be\label{c18}
\|D_p^{m'} \phi_p(D_p)Qs\|_{L^2}  \leq C_{l,m'}
 p^{-l+m}  \|s\|_{L^2}.
\ee
We deduce from (\ref{c17}) and (\ref{c18}) that if $P,Q$ are
differential operators of order $m,m'$  with compact support in
$U_i$, $U_j$ respectively, then for any $l>0$, there exists
$C_l>0$ such that for $p\geq 1$,
\be\label{c19} \|P \phi_p( D_p) Q
s\|_{L^2}\leq C_l p^{-l} \|s\|_{L^2}. \ee
On $U_i\times U_j$, by
using Sobolev inequality and  (\ref{0c3}),
 we get the first inequality of \eqref{1c3}.

 By the finite propagation speed of solutions of hyperbolic equations
\cite{CGT}, \cite{Chernoff}, \cite[\S 7.8]{CP}, \cite[\S 4.4]{T1},
$F(D_p)(x, x')$ only depends on the restriction of $D_p$ to
$B^X(x,\var)$, and is zero if $d(x, x') \geq \var$.
Thus we get the second  inequality of \eqref{1c3}.
The proof of Proposition \ref{0t3.0} is complete.
\end{proof}

From Proposition \ref{0t3.0} and the finite propagation speed as above,
 we know
that the asymptotic of $P_p(x,x')$ as $p\to \infty$ is localized on
a neighborhood of $x$.

To compare the coefficients of the expansion of  $P_p(x,x')$ with the heat kernel  expansion of $\exp(-\frac{u}{p}D^2_p)$ in
Theorem \ref{t0.2}, we will use again the finite propagation speed  to localize the problem.
\begin{defn} For $u>0, a\in \bC$, set
\begin{align}  \label{c4}
&G_u(a)= \int_{-\infty}^{+\infty} e ^{i v a} \exp(-\frac{v^2}{2})
f(\sqrt{u} v) \frac{dv}{\sqrt{2\pi}},\\
&H_u(a) =\int_{-\infty}^{+\infty} e ^{i v a}\exp(-\frac{v^2}{2u})
(1-f(v)) \frac{dv}{ \sqrt{2\pi u}} .\nonumber
\end{align}
\end{defn}

The functions $G_u(a), H_u(a)$ are even holomorphic functions.
The restrictions of $G_u, H_u$
 to $\bR$ lie  in the Schwartz space $\mathcal{S} (\bR)$.
Clearly,
\be\label{c5}
G_{\frac{u}{p}}(\sqrt{\frac{u}{p}} D_p)
+ H_{\frac{u}{p}}(D_p)
= \exp(-\frac{u}{2p} D_p^2).
\ee
Let $G_{\frac{u}{p}}(\sqrt{\frac{u}{p}} D_p)(x,x')$,
$H_{\frac{u}{p}}( D_p)(x,x')$ ($x,x'\in X$)
be the smooth kernels associated to
 $G_{\frac{u}{p}}(\sqrt{\frac{u}{p}} D_p)$,
 $H_{\frac{u}{p}}(D_p)$ calculated  with respect to
the volume form  $dv_X (x')$.

\begin{prop} \label{t3.2} For any $m\in \bN$, $u_0>0, \var>0$,
 there exists $C>0$ such that for any $x,x'\in X$, $p\in \bN$, $u\geq u_0$,
\begin{eqnarray}\label{c6}
\Big |H_{\frac{u}{p}}(D_p)(x,x')\Big|_{\cC^m}
 \leq C p^{2m+2n+2} \exp(-\frac{\var^2 p}{8u}).
\end{eqnarray}
\end{prop}
\begin{proof}
By (\ref{c4}), for any $m\in \bN$ there exists $C_{m}>0$ (which
depends on $\var$) such that
\be\label{c9} \sup_{a\in \bR} |a|^m
|H_u(a) | \leq C_{m} \exp (-\frac{\var^2}{8u}). \ee

As (\ref{c19}), we deduce from (\ref{c17}) and (\ref{c9}) that if
$P,Q$ are differential operators  of order $m,m'$  with compact
support in $U_i$, $U_j$ respectively, then  there exists $ C>0$
such that for $p\geq 1,u\geq u_0$,
\be\label{0c19} \|P
H_{\frac{u}{p}}(D_p) Q s\|_{L^2}\leq C p^{m+m'} \exp
(-\frac{\var^2 p}{8u}) \|s\|_{L^2}. \ee On $U_i\times U_j$, by
using Sobolev inequality, we
 get our Proposition \ref{t3.2}.
\end{proof}

 Using (\ref{c4})
and finite propagation speed \cite[\S 7.8]{CP}, \cite[\S 4.4]{T1},
it is clear that for $x,x'\in X$,
$G_{\frac{u}{p}}(\sqrt{\frac{u}{p}} D_p)(x, x')$ only depends
on the restriction of $D_p$ to
$B^X(x,\var)$, and is zero if $d(x, x') \geq \var$.

\subsection{Rescaling and a Taylor expansion of the operator $D_p$}\label{s3.2}

Now we fix $x_0\in X$.
We identify $L_Z, E_Z$ and $(E_p)_Z$
for $Z\in B^{T_{x_0}X}(0,\var)$ to $L_{x_0}, E_{x_0}$ and $(E_p)_{x_0}$
by parallel transport with respect to the connections
$\nabla ^L, \nabla ^E$ and
 $\nabla^{E_p}$ along the curve $\gamma_Z :[0,1]\ni u \to \exp^X_{x_0} (uZ)$.
 Let $\{e_i\}_i$ be an oriented orthonormal
basis of $T_{x_0}X$.
We also denote by $\{e ^i\}_i$ the dual basis of $\{e_i\}$.
 Let $\wi{e}_i (Z)$  be the parallel
transport of ${e}_i$ with respect to $\nabla^{TX}$
 along the above curve.
%Using $s_L$, we get the isometric $(E_p)_{x_0} \simeq
%(\Lambda (T^{*(0,1)}X)\otimes  E)_{x_0}=\bE_{x_0}$.

Now, for $\var >0$ small enough,  we will extend the geometric objects on
$B^{T_{x_0}X}(0,\var)$ to $\bR^{2n} \simeq T_{x_0}X$
(here we identify $(Z_1,\cdots, Z_{2n})
\in \bR^{2n}$ to $\sum_i Z_i e_i\in T_{x_0}X$) such that
$D_p$ is the restriction of a spin$^c$ Dirac operator on $\bR^{2n}$
associated to a Hermitian line bundle with positive curvature.
In this way, we can replace $X$ by  $\bR^{2n}$.

First of all, we denote $L_0$, $E_0$ the trivial bundles $L_{x_0},
E_{x_0}$ on $X_0= \bR^{2n}$. And we still denote by $\nabla ^L,
\nabla ^E$, $h^L$ etc. the connections and metrics on   $L_0$,
$E_0$ on $B^{T_{x_0}X}(0,4\var)$ induced by the above
identification. Then $h^L$, $h^E$ is  identified with the constant
metrics $h^{L_0}=h^{L_{x_0}}$,  $h^{E_0}=h^{E_{x_0}}$. Let $\mR =
\sum_i Z_i e_i=Z$ be the radial vector field on $\bR^{2n}$.

Let $\rho: \bR\to [0,1]$ be a smooth even function such that
\begin{align}\label{1c14}
\rho (v)=1  \  \  {\rm if} \  \  |v|<2;
\quad \rho (v)=0 \   \   {\rm if} \  |v|>4.
\end{align}
 Let $\varphi_\var : \bR^{2n} \to \bR^{2n}$ is the map defined by
$\varphi_\var(Z)= \rho(|Z|/\var) Z$.
Let $g^{TX_0}(Z)= g^{TX}(\varphi_\var(Z))$, $J_0(Z)=  J(\varphi_\var(Z))$
be the metric and almost complex structure on $X_0$.
 Let $\nabla ^{E_0}= \varphi_\var ^* \nabla ^{E}$, then  $\nabla ^{E_0}$
is the extension of $ \nabla ^{E}$ on $B^{T_{x_0}X}(0,\var)$.
Let $\nabla ^{L_0}$ be the Hermitian connection on $(L_0, h^{L_0})$ defined by
\begin{align}\label{1c15}
&\nabla ^{L_0}|_Z = \varphi_\var ^* \nabla ^{L} +\frac{1}{2}
(1-\rho ^2(|Z|/\var) )  R^{L}_{x_0} (\mR,\cdot).
\end{align}
Then we calculate easily that its curvature $R^{L_0}= (\nabla ^{L_0})^2$ is
\begin{multline}\label{1c16}
R^{L_0}(Z) = \varphi_\var ^* R^L + \frac{1}{2}d \Big((1-\rho ^2(|Z|/\var) )  R^{L}_{x_0} (\mR,\cdot)\Big)\\
= \Big(1-\rho ^2(|Z|/\var)\Big)  R^{L}_{x_0}
+ \rho ^2(|Z|/\var) R^L_{\varphi_\var(Z)}\\
- (\rho \rho')(|Z|/\var)\sum_i \frac{Z_i e ^i}{\var |Z|}\wedge \left[
R^{L}_{x_0} (\mR,\cdot)- R^{L}_{\varphi_\var(Z)} (\mR,\cdot)\right].
\end{multline}
Thus $R^{L_0}$ is positive in the sense of (\ref{b1}) for $\var$ small enough,
 and  the corresponding constant $\mu_0$ for  $R^{L_0}$ is bigger
than $\frac{4}{5}\mu_0$. From now on, we fix  $\var$ as above.

Let $ T^{*(0,1)}X_0$ be the anti-holomorphic cotangent bundle of $(X_0,J_0)$.
Since $J_0(Z)=  J(\varphi_\var(Z))$,
$ T_{Z,J_0}^{*(0,1)}X_0$ is naturally identified with
$T_{\varphi_\var(Z),J}^{*(0,1)}X_0$ (obviously, here the second
subscript indicates the almost complex structure with respect to which the splitting is done).
Let $\nabla ^{\text{Cliff}_0}$ be the Clifford connection on
$\Lambda ( T^{*(0,1)}X_0)$ induced by the Levi-Civita connection
$\nabla ^{TX_0}$
on $(X_0, g^{TX_0})$. Let $R^{E_0}, R^{TX_0}$, $R^{\text{Cliff}_0}$ be the corresponding curvatures on
$E_0,TX_0$ and $\Lambda ( T^{*(0,1)}X_0)$.

 We identify $\Lambda ( T^{*(0,1)}X_0)_Z$ with
$\Lambda ( T^{*(0,1)}_{x_0}X)$ by identifying first
$\Lambda ( T^{*(0,1)}X_0)_Z$
with $\Lambda (T^{*(0,1)}_{\varphi_\var(Z),J}X_0)$, which in turn is
identified with $\Lambda ( T^{*(0,1)}_{x_0}X)$ by using parallel transport
along $u\to u \varphi_\var(Z)$ with respect to
$\nabla ^{\text{Cliff}_0}$. We also trivialize
$\Lambda ( T^{*(0,1)}X_0)$  in this way. Let $S_L$ be an unit vector of $L_{x_0}$. Using $S_L$ and  the above discussion, we
get an isometry $E_{0,p}:=\Lambda ( T^{*(0,1)}X_0)\otimes E_0\otimes L_0^p$
$\simeq (\Lambda ( T^{*(0,1)}X)\otimes E)_{x_0}=: \bE_{x_0}$.

Let $D_p^{X_0}$ (resp. $\nabla^{E_{0,p}}$) be the Dirac operator on $X_0$
(resp. the connection on $E_{0,p}$) associated to the
above data  by the construction in Section \ref{s2}. By the
argument in \cite[p. 656-657]{MM}, we know that Theorem \ref{t2.1}
still holds for $D_p^{X_0}$. In particular, there exists $C>0$
such that
\begin{align}\label{1c17}
&\spec (D_p^{X_0})^2  \subset \{0\}\cup [\frac{8}{5}p\mu_0 -C,+\infty[.
\end{align}

Let $P_p^{0}$ be the orthogonal projection from
$\Omega^{0,\bullet}(X_0,L^p_0\otimes E_0) \simeq \cC^\infty (X_0,\bE_{x_0})$
on $\Ker D_p^{X_0}$, and let  $P_p^{0}(x,x')$ be the smooth kernel of $P_p^{0}$
with respect to the volume form $dv_{X_0}(x')$.
\begin{prop} \label{p3.2} For any $l,m\in \bN$, there exists $C_{l,m}>0$
 such that for $x,x' \in B^{T_{x_0}X}(0,\var)$,
\begin{align}\label{1c19}
\Big |(P_p^{0}- P_p)(x,x')\Big |_{\cC^m}\leq C_{l,m} p^{-l}.
\end{align}
\end{prop}
\begin{proof} Using (\ref{1c2}) and (\ref{1c17}), we know that
$P_p^{0}-F(D_p)$ verifies also (\ref{1c3})
for $x,x' \in B^{T_{x_0}X}(0,\var)$, thus we get (\ref{1c19}).
\end{proof}

To be complete, we prove the following result in
\cite[Proposition 1.28]{BeGeV}.

\begin{lemma} \label{l3.2} The Taylor expansion of $\wi{e}_i(Z)$
with respect to the basis $\{e_i\}$ to order $r$ is a polynomial of
 the Taylor expansion of the coefficients of $R^{TX}$
to order $r-2$. Moreover we have
\be\label{0c21}
\wi{e}_i(Z)  = e_i - \frac{1}{6} \sum_{j}
\left \langle R^{TX}_{x_0} (\mR,e_i) \mR, e_j\right \rangle  e_j
+ \sum_{|\alpha|\geq 3}
\Big (\frac{\partial^\alpha}{ \partial Z^\alpha}\wi{e}_i\Big ) (0)\frac{Z^\alpha}{\alpha !}.
\ee
\end{lemma}
\begin{proof}
Let $\Gamma ^{TX}$ be the connection form of  $\nabla^{TX}$
 with respect to
the frame $\{\wi{e}_i\}$ of $TX$.
Let $\partial_i=\nabla _{e_i}$  be the partial derivatives along $e_i$.
By the definition of our fixed frame, we have
$i_{\mR}\Gamma ^{TX} =0$.
As in  \cite[(1.12)]{BeGeV},
\begin{align}\label{0c24}
\mL_\mR  \Gamma ^{TX} = [i_{\mR}, d] \Gamma ^{TX}
= i_{\mR} (d \Gamma ^{TX} + \Gamma ^{TX} \wedge \Gamma ^{TX})
= i_{\mR} R^{TX}.
\end{align}

Let $\Theta (Z) = (\theta _j^i (Z))_{i,j=1}^{2n}$
be the $2 n \times 2n$-matrix such that
\be\label{0c25}
e_i = \sum_j \theta ^j_i(Z) \wi{e}_j (Z), \quad
\wi{e}_j (Z)= (\Theta (Z)^ {-1})_j^k e_k.
\ee
Set $\theta ^j (Z) = \sum_i \theta ^j_i(Z) e^i$ and
\begin{align}\label{0c26}
&\theta = \sum_j e^j \otimes e_j
= \sum_j \theta ^j \wi{e}_j \in T^*X\otimes TX.
\end{align}

As $\nabla ^{TX}$ is torsion free, $ \nabla ^{TX}\theta=0$,
thus the $\bR ^{2n}$-valued one-form $\theta= (\theta ^j(Z))$
 satisfies the structure equation,
\be\label{0c27}
d \theta + \Gamma ^{TX} \wedge \theta =0.
\ee
Observe first that (cf. \cite[Proposition 1.27]{BeGeV})
\begin{align}\label{0c28}
&\mR= \sum_j Z_j \wi{e}_j (Z),\quad
i_\mR \theta = \sum_j Z_j e_j = \mR.
\end{align}

Substituting  (\ref{0c28}) and
$(\mL _{\mR} -1) \mR =0$, into the identity $i_\mR (d \theta +\Gamma ^{TX} \wedge \theta)=0$, we obtain
\begin{align}\label{0c29}
(\mL _{\mR} -1) \mL _{\mR} \theta =
(\mL _{\mR} -1) ( d\mR + \Gamma ^{TX} \mR)
= (\mL _{\mR} \Gamma ^{TX}) \mR
= (i_\mR R^{TX})\mR.
\end{align}
Using (\ref{0c28}) once more gives
\begin{align}\label{0c30}
i_{e_j} (\mL _{\mR} -1) \mL _{\mR} \theta ^i (Z)=
\left \langle R^{TX} ( \mR,e_j) \mR, \wi{e}_i\right \rangle (Z) .
\end{align}
Thus
\begin{align}\label{0c31}
\sum_{|\alpha| \geq 1} ( |\alpha|^2 + |\alpha|)
(\partial ^\alpha\theta ^i_j)(0) \frac{Z^\alpha}{\alpha !}=
\left \langle R^{TX} ( \mR,e_j) \mR, \wi{e}_i\right \rangle (Z) .
\end{align}

Now by (\ref{0c25})  and $\theta ^i_j (x_0) = \delta_{ij}$,
 (\ref{0c31}) determines the Taylor expansion of $\theta ^i_j(Z)$
to order $m$ in terms of the Taylor expansion of $R^{TX}$ to order $m-2$.
And
\begin{align}\label{0c32}
(\Theta ^{-1}) ^i_j  = \delta_{ij}
- \frac{1}{6}
\left \langle R^{TX}_{x_0} (\mR,e_i) \mR, e_j\right \rangle
 + \cO(|Z|^3).
\end{align}
By (\ref{0c25}), (\ref{0c32}), we get  (\ref{0c21}).
\end{proof}

For  $s \in \cC^{\infty}(\bR^{2n}, \bE_{x_0})$ and $Z\in \bR^{2n}$,
for $t=\frac{1}{\sqrt{p}}$, set
\begin{align}\label{c27}
&(S_{t} s ) (Z) =s (Z/t),\quad    \nabla_{t}=  S_t^{-1}
t\nabla ^{E_{0,p}} S_t, \\
 &{\bf D}_t= S_t^{-1}  t D_p^{X_0} S_t,
\quad   L^t_2= S_t^{-1}  t^2 D_p^{X_0,2} S_t.\nonumber
\end{align}
 Denote by  $\nabla_U$ the ordinary differentiation
 operator on $T_{x_0}X$ in the direction $U$.
If $\alpha = (\alpha_1,\cdots, \alpha_{2n})$ is a multi-index,
set $Z^\alpha = Z_1^{\alpha_1}\cdots Z_{2n}^{\alpha_{2n}}$.
Set
\begin{align}\label{c29}
\mO_0 = \sum_j  c(e_j)
\Big(\nabla_{e_j}+\frac{1}{2}  R^L_{x_0}(Z, e_j)\Big).
\end{align}

\begin{thm}\label{t3.3} There exist $\mB_{i,r}$
 {\rm (}resp. $\mA_{i,r}$, resp. $\mC_{i,r}${\rm )}
{\rm (}$r\in \bN, i\in \{1,\cdots, 2n\}${\rm )} homogeneous polynomials in $Z$
of degree $r$ with coefficients polynomials in $R^{TX}$,
$R^{\det}$, $R^E$ {\rm (}resp. $R^{TX}$, resp.  $R^L$, $R^{TX}${\rm )} and
their derivatives at $x_0$ to order $r-1$ {\rm (}resp. $r-2$, resp.
$r-1$, $r-2${\rm )}
 such that if we denote by
\begin{align}\label{0c35}
\mO_r = \sum_{i=1}^{2n} c(e_i) \Big( \mA_{i,r}\nabla_{e_i}
+ \mB_{i,r-1}+ \mC_{i,r+1}\Big),
\end{align}
 then
\begin{align}\label{c30}
{\bf D}_t= \mO_0
+ \sum_{r=1}^m t^r \mO_r + \cO(t^{m+1}).
\end{align}
Moreover, there exists $m'\in \bN$ such that for any $k\in \bN$, $t\leq 1$,
$|tZ|\leq \var$,
the derivatives of order $\leq k$ of the coefficients of the operator
 $\cO(t^{m+1})$ are dominated by $C t^{m+1} (1+|Z|)^{m'}$.
\end{thm}
\begin{proof}  By the definition of $\nabla^{{\rm Cliff}}$, $\wi{e}_j$,
for $Z\in \bR^{2n}$, $|Z|\leq \var$,
\be\label{c31}
[\nabla^{\text{Cliff}}_Z, c(\wi{e}_j)(Z)] = c(\nabla ^{TX}_Z\wi{e}_j)(Z)=0.
\ee
Thus  we know that under our trivialization,
for $Z\in \bR^{2n}$, $|Z|\leq \var$,
\be\label{c32}
 c(\wi{e}_j)(Z)=c(e_j).
\ee

We identify $(\det (T^{(1,0)}X))_Z$ for $Z\in B^{T_{x_0}X}(0,\var)$ to
$(\det (T^{(1,0)}X))_{x_0}$
by parallel transport with respect to the connection
$\nabla^{\det (T^{(1,0)}X)}$ along the curve $\gamma_Z$.
Let $\Gamma ^E$, $\Gamma ^{\det}$ and  $\Gamma ^L$
be the connection forms of $\nabla^E$,
 $\nabla^{\det (T^{(1,0)}X)}$ and $\nabla^L$
 with respect to any fixed frames for
$E$, $\det (T^{(1,0)}X)$ and $L$
 which are parallel along the curve $\gamma_Z$
under our trivialization on $B^{T_{x_0}X}(0,\var)$. Then the
corresponding connection form of $\Lambda (T^{*(0,1)}X)$ is
\be\label{0c36} \Gamma^{\text{Cliff}} = \frac{1}{4} \left \langle
\Gamma ^{TX} \wi{e}_k, \wi{e}_l  \right \rangle  c(\wi{e}_k)
c(\wi{e}_l) + \frac{1}{2} \Gamma ^{\det}. \ee

Now for $\Gamma ^\bullet= \Gamma ^E, \Gamma ^L$
or $\Gamma ^{\det}$ and
 $R^\bullet= R^E, R^L$ or $R^{\det}$ respectively,
by the definition of our fixed frame, we have as in (\ref{0c24})
\begin{align}\label{c33}
&i_{\mR}\Gamma ^\bullet =0,
&\mL_\mR  \Gamma ^\bullet = [i_{\mR}, d] \Gamma ^\bullet
= i_{\mR} (d \Gamma ^\bullet + \Gamma ^\bullet \wedge \Gamma ^\bullet)
= i_{\mR} R^\bullet.
\end{align}
Expanding  the Taylor's series of both sides of (\ref{c33}) at $Z=0$, we obtain
\begin{align}\label{c34}
\sum_{\alpha} (|\alpha|+1)  (\partial^\alpha \Gamma ^\bullet )_{x_0} (e_j)
\frac{Z^\alpha}{\alpha !}
= \sum_{\alpha} (\partial^\alpha R^\bullet )_{x_0}(\mR, e_j)
  \frac{Z^\alpha}{\alpha !}.
\end{align}
By equating coefficients of $Z^\alpha$ on both sides, we see from
this formula
\begin{align}\label{c35}
\sum_{|\alpha|=r}  (\partial^\alpha \Gamma ^\bullet )_{x_0} (e_j)
 \frac{Z^\alpha}{\alpha !}
=\frac{1}{r+1} \sum_{|\alpha|=r-1}
(\partial^\alpha R^\bullet )_{x_0}(\mR, e_j)
  \frac{Z^\alpha}{\alpha !}
\end{align}
Especially,
\be\label{c36}
\partial_i\Gamma ^\bullet_{x_0}(e_j) = \frac{1}{2} R^\bullet_{x_0}(e_i, e_j).
\ee
Furthermore, it follows that  the Taylor coefficients of
$\Gamma ^\bullet (e_j) (Z)$ at $x_0$
to order $r$ are determined by those of $R^\bullet$ to order $r-1$.

By (\ref{c29}), (\ref{c32}), for $t= 1/\sqrt{p}$,
 for $|Z|\leq \sqrt{p} \var$, then
\begin{align}\label{c37}
&\nabla_t |_Z  = \nabla + (t \Gamma ^{\rm Cliff}
+ t \Gamma ^E  + \frac{1}{t} \Gamma ^L)  (tZ),\\
&{\bf D}_t
= \sum_{j=1}^{2n} c(e_j) \nabla_{t,\wi{e}_j(tZ)} |_Z.\nonumber
\end{align}
By Lemma \ref{l3.2},
(\ref{c35}) and (\ref{c37}),  we get our Theorem.
\end{proof}

\subsection{Uniform estimate on the heat kernel and the Bergman kernel}
\label{s3.3}

Recall that the operators $L^t_2$, $\nabla_{t}$ were defined in (\ref{c27}).
We also denote by $\left \langle \ ,\ \right\rangle_{0,L^2}$ and $\|\ \|_{0,L^2}$
the scalar product and the $L^2$ norm on $\cC^\infty (X_0, \bE_{x_0})$
induced by $g^{TX_0}, h^{E_0}$ as in (\ref{b3}).

Let $dv_{TX}$ be the Riemannian volume form on
$(T_{x_0}X, g^{T_{x_0}X})$.
Let $\kappa (Z)$ be the smooth positive function defined by the equation
\be\label{c22}
dv_{X_0}(Z) = \kappa (Z) dv_{TX}(Z),
\ee
with $k(0)=1$.
For $s\in \cC^{\infty}( T_{x_0}X, \bE_{x_0}) $, set
\begin{align}\label{u0}
&\|s\|_{t,0}^2 = \int_{\bR^{2n}} |s(Z)|^2_{h^{\Lambda ( T^{*(0,1)}X_0)\otimes E_0}(tZ)}dv_{X_0}(tZ)= t^{-2n}\|S_t s\|_{0,L^2}^2,\\
&\| s \|_{t,m}^2 = \sum_{l=0}^m \sum_{i_1,\cdots, i_l=1}^{2n}
\| \nabla_{t,e_{i_1}} \cdots \nabla_{t,e_{i_l}} s\|_{t,0}^2. \nonumber
\end{align}
We denote by $\left \langle s ', s \right\rangle_{t,0}$
 the inner product on $\cC^\infty (X_0, \bE_{x_0})$
 corresponding to $\|\quad\|^2_{t,0}$.
Let $H^m_t$ be the Sobolev space of order $m$ with norm $\|\quad\|_{t,m}$.
Let $H^{-1}_t$ be the Sobolev space of order $-1$ and let $\|\quad\|_{t,-1}$ be the norm on  $H^{-1}_t$ defined by
$\|s\|_{t,-1} = \sup_{0\neq s'\in  H^1_t }
|\left \langle s,s'\right \rangle_{t,0}|/\|s'\|_{t,1}$.
If $A\in \cL (H^{m}, H^{m'})$ $(m,m' \in \bZ)$,
 we denote
by  $\| A \|^{m,m'}_t$
 the norm of $A$ with respect to the norms
$\|\quad  \|_{t,m}$ and $\| \quad \|_{t,m'}$.

Then $L^t_2$ is a formally self adjoint elliptic operator with respect to
 $\|\quad\|^2_{t,0}$, and is a smooth family of
operators with parameter $x_0\in X$.

%Let $\Delta^{TX} = \sum_i \frac{\partial ^2}{\partial Z_i^2}$
%be the standard Euclidean Laplacian on $T_{x_0}X$ with
%respect to the metric $g^{T_{x_0}X}$.

\begin{thm}  \label{tu1} There exist constants $ C_1, C_2, C_3>0$
such that for $t\in ]0,1]$ and any
$s,s'\in \cC^{\infty}_0(\bR^{2n}, \bE_{x_0})$,
\begin{align}\label{u1}
& \left \langle L^t_2 s,s\right \rangle_{t,0}
\geq C_1\|s\|_{t,1}^2 -C_2 \|s\|_{t,0}^2  , \\
& |\left \langle L^t_2 s, s'\right \rangle_{t,0}|
 \leq C_3 \|s\|_{t,1}\|s'\|_{t,1}.\nonumber
\end{align}
\end{thm}

\begin{proof}
Now from (\ref{Lich}),
\begin{align}\label{1u1}
 \left \langle D_p^{X_0,2}s,s\right \rangle_{0,L^2} =
\|\nabla ^{E_{0,p}}s\|_{0,L^2}^2
+ \left \langle \left (-2p \om_d
-p\tau+\tfrac{1}{4}r^X+ \mathbf{c}(R)\right )s,s\right \rangle_{0,L^2}.
\end{align}
Thus from \eqref{c27}, (\ref{u0}), and (\ref{1u1}),
\begin{align}\label{1u2}
\left \langle  L^t_2 s,s\right \rangle_{t,0} = \|\nabla_ts\|_{t,0}^2
+ \left \langle\left  (-2 S_t^{-1}\om_d
-S_t^{-1}\tau + \tfrac{t^2}{4}S_t^{-1}r^X+t^2 S_t^{-1}\mathbf{c}(R)\right )s,s\right \rangle_{t,0}.
\end{align}
From (\ref{1u2}), we get (\ref{u1}).
\end{proof}

Let $\delta$ be the counterclockwise oriented circle in $\bC$
of center $0$ and radius $\mu_0/4$,
and let $\Delta$ be the oriented path in $\bC$ which goes parallel to the real
axis from $+\infty +i$ to $\frac{\mu_0}{2}+i$ then parallel to the imaginary
 axis to $\frac{\mu_0}{2}-i$ and the parallel to the real
axis to $+\infty -i$.
By (\ref{1c17}), (\ref{c27}), for $t$ small enough,
\begin{align}\label{1u3}
\spec \, L^t_2\subset \{0\}\cup [\mu_0,+ \infty[.
\end{align}
Thus $(\lambda- L^t_2)^{-1}$ exists for $\lambda \in \delta\cup \Delta$.

\begin{thm}\label{tu4}  There exists $C>0$ such that for  $t\in ]0,1]$,
$\lambda \in  \delta\cup \Delta$, and $x_0\in X$,
\begin{align}\label{ue2}
& \| (\lambda- L^t_2)^{-1}\|^{0,0}_t \leq C,\\
& \| (\lambda- L^t_2)^{-1}\|^{-1,1}_t
\leq C (1+|\lambda |^2).\nonumber
\end{align}
\end{thm}
\begin{proof} The  first inequality of (\ref{ue2}) is from (\ref{1u3}).
Now, by (\ref{u1}), for $\lambda_0\in\bR$, $\lambda_0\leq -2C_2$,
 $(\lambda_0- L^t_2)^{-1}$ exists, and we have
$\|(\lambda_0- L^t_2)^{-1}\|^{-1,1}_t \leq \frac{1}{C_1}$. Now,
\be\label{ue7}
(\lambda- L^t_2)^{-1}= (\lambda_0- L^t_2)^{-1}
- (\lambda-\lambda_0) (\lambda- L^t_2)^{-1}(\lambda_0- L^t_2)^{-1}.
\ee
Thus for  $\lambda\in \delta\cup \Delta$, from (\ref{ue7}), we get
\be\label{ue8}
\|(\lambda-L^t_2)^{-1}\|^{-1,0}_t  \leq
\frac{1}{C_1} \left(1+\frac{4}{\mu_0}|\lambda-\lambda_0|\right).
\ee
Now we  change the last two factors in (\ref{ue7}), and apply (\ref{ue8}),
we get
\begin{align}\label{ue9}
\|(\lambda-L^t_2)^{-1}\|^{-1,1}_t & \leq  \frac{1}{C_1}+
 \frac{|\lambda-\lambda_0|}{{C_1}^2}
\left(1+\frac{4}{\mu_0}|\lambda-\lambda_0|\right)\\
&\leq C(1+|\lambda|^2).\nonumber
\end{align}
The proof of our Theorem is complete.
\end{proof}

\begin{prop} \label{tu5} Take $m \in \bN^*$. There exists $C_m>0$ such that
 for
  $t\in ]0,1]$,  $Q_1, \cdots$, $Q_m \in$ $\{ \nabla_{t,e_i}, Z_i\}_{i=1}^{2n}$
and  $s,s'\in \cC^{\infty}_{0}(\bR^{2n}, \bE_{x_0})$,
\be\label{ue11}
\left |\left \langle [Q_1,
[Q_2,\ldots, [Q_m,  L^t_2]] \ldots ]s,
s'\right \rangle_{t,0} \right |
\leq C_m  \|s\|_{t,1} \|s'\|_{t,1}.
\ee
\end{prop}
\begin{proof} Set $g_{ij}(Z)=g^{TX_0} (e_i,e_j)(Z)$.
 Let $(g^{ij}(Z))$ be the inverse of the matrix $(g_{ij}(Z))$.
 Let $\nabla ^{TX_0}_{e_i}e_j= \Gamma_{ij}^k(Z)e_k$, then by (\ref{Lich}),
\begin{multline}\label{1ue1}
 L^t_2(Z) = - g^{ij}(tZ) (\nabla_{t,e_i}\nabla_{t,e_j}-t \Gamma_{ij}^k(tZ) \nabla_{t,e_i}) \\
-2  \om_d(tZ)+ \tau(tZ)+ t^2 (\tfrac{1}{4}r^X+ \mathbf{c}(R))(tZ).
\end{multline}
Note that $[\nabla_{t,e_i}, Z_j]=\delta_{ij}$.
Thus by (\ref{1ue1}), we know that $[Z_j, L^t_2]$ verifies (\ref{ue11}).

Note that by  (\ref{c27}),
\begin{align}\label{1ue2}
[\nabla_{t,e_i},\nabla_{t,e_j}]= \left(R^{L_0}
+ t^2  R^{\text{Cliff}_0} +   t^2 R^{E_0}\right)(tZ)(e_i,e_j).
\end{align}
Thus from (\ref{1ue1}) and (\ref{1ue2}), we know that
$[\nabla_{t,e_k},  L^t_2]$
has the same structure as  $L^t_2$ for $t \in ]0,1]$,
i.e. $[\nabla_{t,e_k},  L^t_2]$ has the type as
\begin{align}\label{1ue3}
\sum_{ij} a_{ij}(t,tZ) \nabla_{t,e_i}\nabla_{t,e_j}
+\sum_{i}b_{i}(t,tZ) \nabla_{t,e_i} + c(t,tZ),
\end{align}
and $a_{ij}(t,Z),b_{i}(t,Z),  c(t,Z)$ and their derivatives on $Z$
are uniformly bounded for $Z\in \bR ^{2n}, t\in [0,1]$;
moreover, they are polynomial in $t$.

Let $(\nabla_{t,e_i})^*$ be the adjoint of $\nabla_{t,e_i}$ with respect to
$\left \langle \ , \ \right \rangle_{t,0}$, then by (\ref{u0}),
\begin{align}\label{1ue4}
(\nabla_{t,e_i})^* =- \nabla_{t,e_i}
- t (k^ {-1}\nabla_{e_i}k)(tZ),
\end{align}
the last term of (\ref{1ue4}) and its derivatives in  $Z$ are uniformly bounded in $Z\in \bR ^{2n}, t\in [0,1]$.

By (\ref{1ue3}) and (\ref{1ue4}), (\ref{ue11}) is verified for
$m=1$.

By iteration, we know that
$[Q_1,[Q_2,\ldots, [Q_m,  L^t_2]] \ldots ]$
has the same structure (\ref{1ue3})  as $L^t_2$.
By (\ref{1ue4}),
we get Proposition \ref{tu5}.
\end{proof}

\begin{thm}\label{tu6} For any $t\in ]0,1]$, $\lambda \in \delta\cup \Delta$,
 $m \in \bN$, the resolvent $(\lambda-L^t_2)^{-1}$ maps $H^m_t$
into $H^{m+1}_t$. Moreover for any $\alpha\in \bZ^{2n}$,
 there exist $N \in \bN$, $C_{\alpha, m}>0$
such that for  $t\in ]0,1]$, $\lambda \in \delta\cup \Delta$,
 $s\in \cC^\infty_0 (X_0,\bE_{x_0})$,
\be\label{ue12}
\|Z^\alpha (\lambda-L^t_2)^{-1} s\|_{t,m+1}
\leq C_{\alpha, m}  (1+|\lambda|^2)^N \sum_{\alpha' \leq \alpha}
\|Z^{\alpha'} s\|_{t,m} .
\ee
\end{thm}
\begin{proof} For $Q_1, \cdots, Q_m\in \{\nabla_{t,e_i}\}_{i=1}^{2n}$,
$Q_{m+1},\cdots, Q_{m+|\alpha|}\in\{Z_i\}_{i=1}^{2n}$,  We can express
$Q_1\cdots$ $Q_{m+|\alpha|}(\lambda-L^t_2)^{-1}$ as a linear combination of operators of the type
\be\label{ue13}
[Q_1, [ Q_2,\ldots [Q_{m'},(\lambda-L^t_2)^{-1}]]\ldots ]Q_{m'+1}
\cdots Q_{m+|\alpha|} \quad m'\leq m+|\alpha|.
\ee
Let $\cR_{t}$ be the family operators
$\cR_{t} = \{ [ Q_{j_1}, [Q_{j_2},\ldots [Q_{j_l},L^t_2]]\ldots ] \}.
$
Clearly, any commutator
$[Q_1, [ Q_2,\ldots [Q_{m'},(\lambda-L^t_2)^{-1}]]\ldots ]$
is a linear combination of operators of the form
\be\label{ue14}
(\lambda-L^t_2)^{-1}R_1(\lambda-L^t_2)^{-1}R_2
\cdots R_{m'}(\lambda-L^t_2)^{-1}
\ee
with  $R_1, \cdots, R_{m'} \in \cR_{t}$.

By Proposition \ref{tu5}, the norm $\|\quad\|_t^{1,-1}$ of the
 operators $R_j\in \cR_{t}$ is
uniformly bound by $C$.
By Theorem \ref{tu4}, we find that there exist $C>0$, $N \in \bN$
such that the norm $\|\quad\|_t^{0,1}$
of operators (\ref{ue14}) is dominated by
$C  (1+|\lambda|^2)^N$.
\end{proof}

Let $e ^{-uL^t_2}(Z,Z')$, $(L^t_2e ^{-uL^t_2})(Z,Z')$
be the smooth kernels of the operators  $e ^{-uL^t_2}$, $L^t_2e ^{-uL^t_2}$
with respect to $dv_{TX}(Z')$.
Note that $L^t_2$ are families of differential operators with coefficients in
$\End ( \bE_{x_0}) =\End (\Lambda (T^{*(0,1)}X)\otimes E)_{x_0}$.
Let $\pi : TX\times_{X} TX \to X$ be the natural projection from the
fiberwise product of $TX$ on $X$.
Then we can view  $e ^{-uL^t_2}(Z,Z')$, $(L^t_2e ^{-uL^t_2})(Z,Z')$
as smooth sections
of $\pi ^* (\End (\Lambda (T^{*(0,1)}X)\otimes E))$ on  $TX\times_{X} TX$.
 Let $\nabla ^{\End (\bE)}$ be the connection on
$\End (\Lambda (T^{*(0,1)}X)\otimes E)$ induced by
$\nabla ^{\text{Cliff}}$ and $\nabla ^E$. And $\nabla ^{\End (\bE)}$
induces naturally a $\cC^m$-norm for the parameter $x_0\in X$.

\begin{thm}\label{tue8} There exist $C'' >0$ such that
 for any $m,m', r\in \bN$, $u_0>0$, there exist $C>0$, $N\in \bN$
 such that for $t \in ]0,1]$, $u\geq u_0$, $Z,Z'\in T_{x_0}X$,
 \begin{align}\label{ue15}
&\sup_{|\alpha|,|\alpha'|\leq m}
\Big |\frac{\partial^{|\alpha|+|\alpha'|}}
{{\partial Z}^{\alpha} {\partial Z'}^{\alpha'}}
\frac{\partial^{r}}{\partial t^{r}}
e ^{-uL^t_2} \left (Z, Z'\right )\Big |_{\cC^{m'}(X)} \\
&\hspace{20mm} \qquad   \leq C (1+|Z|+|Z'|)^N
\exp (\frac{1}{2}\mu_0 u- \frac{2C''}{u} |Z-Z'|^2), \nonumber\\
&\sup_{|\alpha|,|\alpha'|\leq m}
\Big |\frac{\partial^{|\alpha|+|\alpha'|}}
{{\partial Z}^{\alpha} {\partial Z'}^{\alpha'}}
\frac{\partial^{r}}{\partial t^{r}}
(L^t_2e ^{-uL^t_2})  \left (Z, Z'\right )\Big |_{\cC^{m'}(X)} \nonumber \\
&\hspace{20mm} \qquad   \leq C (1+|Z|+|Z'|)^N
\exp ( -\frac{1}{4}\mu_0 u- \frac{2C''}{u} |Z-Z'|^2). \nonumber
\end{align}
here $\cC^{m'}(X)$ is the $\cC^{m'}$ norm for the parameter $x_0\in X$.
\end{thm}
\begin{proof} By (\ref{1u3}), for any $k\in \bN^*$,
\begin{align}\label{1ue15}
&e ^{-uL^t_2}= \frac{(-1)^{k-1} (k-1)!}{2\pi i u^{k-1}}
\int_{\delta\cup\Delta} e^{-u\lambda} (\lambda - L^t_2)^{-k} d \lambda,\\
&L^t_2e ^{-uL^t_2}= \frac{(-1)^{k-1} (k-1)!}{2\pi i u^{k-1}}
\int_{\Delta} e^{-u\lambda}\Big[ \lambda (\lambda - L^t_2)^{-k}
-  (\lambda - L^t_2)^{-k+1} \Big]  d \lambda.\nonumber
\end{align}
For $m\in \bN$, let $\mQ ^m$ be the set of operators
  $\{\nabla_{t,e_{i_1}}\cdots  \nabla_{t,e_{i_j}}\}_{j\leq m}$.
From Theorem \ref{tu6}, we deduce that
if $Q\in \mQ ^{m}$, there are $M\in \bN$, $C_m>0$ such that
for any $\lambda \in \delta\cup \Delta$,
\be\label{ue18}
\| Q (\lambda-L^t_2)^{-m}\|_t^{0,0}
 \leq C_m  (1+|\lambda |^2) ^M.
\ee
Next we study $L^{t*}_2$, the formal adjoint of $L^t_2$
with respect to  (\ref{u0}).
 Then $L^{t*}_2$
has  the same structure (\ref{1ue3})  as the operator $L^t_2$,
especially, \be\label{ue19} \| Q (\lambda-L^{t*}_2)^{-m}\|_t^{0,0}
\leq C_m  (1+|\lambda|^2) ^M. \ee After taking the adjoint of
(\ref{ue19}),  we get \be\label{ue20} \|
(\lambda-L^t_2)^{-m}Q\|_t^{0,0} \leq C_m  (1+|\lambda |^2)^M. \ee
From  (\ref{1ue15}), (\ref{ue18}) and (\ref{ue20}), we have, for
$Q,Q' \in \mQ ^{m}$,
\begin{align}\label{ue21}
&\|Q e ^{-uL^t_2}Q' \|^{0,0}_t \leq C_{m} e^{\frac{1}{4}\mu_0 u},\\
&\|Q (L^t_2e ^{-uL^t_2})Q' \|^{0,0}_t \leq C_{m} e^{-\frac{1}{2}\mu_0 u}.\nonumber
\end{align}

Let $|\quad|_m$ be the usual Sobolev norm on $\cC^\infty(\bR^{2n}, \bE_{x_0})$
induced by $h^{\bE_{x_0}}= h^{\Lambda (T^{*(0,1)}_{x_0}X)\otimes E_{x_0}}$
and the volume form $dv_{TX}(Z)$ as in (\ref{u0}).
Observe that by (\ref{c37}), (\ref{u0}), there exists $C>0$ such that
for $s\in \cC^\infty (X_0, \bE_{x_0})$,
$\supp s \subset B^{T_{x_0}X}(0,q)$, $m\geq 0$,
\begin{align}\label{1ue21}
&\frac{1}{C} (1+q)^{-m} \|s\|_{t,m}\leq |s|_m \leq C(1+q)^m \|s\|_{t,m}.
\end{align}
Now (\ref{ue21}), (\ref{1ue21})  together with Sobolev's inequalities implies that
if $Q,Q' \in \mQ ^{m}$,
\begin{align}\label{ue23}
&\sup_{|Z|,|Z'|\leq q}| Q_Z Q'_{Z'}
e ^{-uL^t_2} (Z,Z')  | \leq C(1+q)^{2n+2} \ e^{\frac{1}{4}\mu_0 u},\\
&\sup_{|Z|,|Z'|\leq q}| Q_Z Q'_{Z'}
 (L^t_2e ^{-uL^t_2}) (Z,Z')  | \leq C(1+q)^{2n+2} \ e^{-\frac{1}{2}\mu_0 u}.
\nonumber
\end{align}
Thus by (\ref{c37}), (\ref{ue23}), we derive (\ref{ue15})
with the exponential $e^{\frac{1}{4}\mu_0 u}$,
$e^{-\frac{1}{2}\mu_0 u}$
for the case when $r =m'=0$ and $C''=0$.

To obtain (\ref{ue15}) in general, we proceed as in the proof of \cite[Theorem 11.14]{B95}.
Note that the function $f$ is defined in (\ref{c2}). For $h >1$,
put
\begin{align}  \label{ue24}
&K_{u,h}(a)= \int_{-\infty}^{+\infty} \exp(i v\sqrt{2u} a)
 \exp(-\frac{v^2}{2})
\Big (1-f (\frac{1}{h}\sqrt{2u} v) \Big ) \frac{dv}{\sqrt{2\pi}}.
\end{align}
Then there exist $C',C_1>0$ such that for any $c>0$, $m,m'\in \bN$,
there is $C>0$ such that for $t\in ]0,1], u\geq u_0$, $h>1$,
$a\in \bC, |{\rm Im} (a)|\leq c$,  we have
\begin{align} \label{1ue24}
|a|^m |K_{u,h}^{(m')}(a)| \leq C \exp \Big( C'c^2 u- \frac{C_1}{u} h^2\Big).
\end{align}
For any $c>0$, let  $V_c$ be the images of
$\{\lambda \in \bC, |\Im \lambda |\leq c\}$
by the map $\lambda \to \lambda ^2$. Then
$
V_c=\{\lambda \in \bC, {\rm Re} (\lambda) \geq \frac{1}{4c^2}
{\rm Im} (\lambda)^2 -c^2\}$, and $\delta \cup \Delta\subset V_c$
for $c$ big enough.
Let $\wi{K}_{u,h}$ be the holomorphic function such that
$\wi{K}_{u,h}(a^2)=K_{u,h}(a)$. Then by (\ref{1ue24}), for $\lambda \in V_c$,
\be\label{ue25}
|\lambda|^m |\wi{K}_{u,h}^{(m')}(\lambda)|
\leq C \exp\Big(C'c^2 u - \frac{C_1}{u} h^2\Big).
\ee
Using finite propagation speed of solutions of hyperbolic equations
 and (\ref{ue24}), we find that there exists a fixed constant
(which depends on $\varepsilon$) $c'>0$ such that
\begin{align}\label{1ue25}
\wi{K}_{u,h}(L^t_2)(Z,Z') = e ^{-u L^t_2}(Z,Z')\quad {\rm if} \,\,
|Z-Z'|\geq c' h.
\end{align}
By (\ref{ue25}), we see that given $k\in \bN$, there is a unique
holomorphic function $\wi{K}_{u,h,k} (\lambda)$ defined on
a neighborhood of $V_c$ such that  it verifies the same estimates as $\wi{K}_{u,h}$ in (\ref{ue25}) and
$\wi{K}_{u,h,k} (\lambda) \to 0$ as $\lambda \to +\infty$; moreover
\be\label{ue16}
\wi{K}_{u,h,k}^{(k-1)} (\lambda)/(k-1)! = \wi{K}_{u,h} (\lambda).
\ee
Thus as in (\ref{1ue15}),
\begin{align} \label{1ue26}
&\wi{K}_{u,h}(L^t_2)= \frac{1}{2\pi i} \int_{\delta\cup\Delta}\wi{K}_{u,h,k}
 (\lambda)(\lambda-L^t_2)^{-k} d \lambda,\\
&L^t_2\wi{K}_{u,h}(L^t_2)= \frac{1}{2\pi i} \int_{\Delta}
\wi{K}_{u,h,k} (\lambda)\Big[\lambda (\lambda-L^t_2)^{-k}
- (\lambda-L^t_2)^{-k+1}\Big]  d \lambda. \nonumber
\end{align}

By (\ref{ue18}), (\ref{ue20}) and by proceeding
as in (\ref{ue21})-(\ref{ue23}),
we find that for ${\bf K}(a)= \wi{K}_{u,h}(a)$
or $a\wi{K}_{u,h}(a)$, for  $|Z|, |Z'|\leq q$,
 \begin{align}  \label{ue27}
\sup_{|\alpha|,|\alpha'|\leq m}
\Big |\frac{\partial^{|\alpha|+|\alpha'|}}
{\partial Z^{\alpha} {\partial Z'}^{\alpha'}}
{\bf K}(L^t_2) (Z, Z')\Big |
 \leq C (1+q)^N
\exp \left(C'c^2 u- \frac{C_1}{u} h^2\right).
\end{align}
Setting $h= \frac{1}{c'}|Z-Z'|$ in (\ref{ue27}), we get
for $\alpha, \alpha'$ verified $|\alpha|,|\alpha'|\leq m$,
\begin{align}  \label{1ue27}
\Big |\frac{\partial^{|\alpha|+|\alpha'|}}
{\partial Z^{\alpha} {\partial Z'}^{\alpha'}}
{\bf K}(L^t_2) (Z, Z')\Big |
 \leq C (1+|Z|+|Z'|)^N
\exp \left(C'c^2  u- \frac{C_1}{2{c'}^2u} |Z-Z'|^2\right).
\end{align}
By (\ref{ue15})
with the exponential $e^{\frac{1}{4}\mu_0 u}$, $e^{-\frac{1}{2}\mu_0 u}$
for $r=m'=C''=0$, (\ref{1ue25}), (\ref{1ue27}),
we get (\ref{ue15}) for $r=m'=0$.

To get (\ref{ue15}) for $r \geq 1$, note that from (\ref{1ue15}), for $k\geq 1$
 \begin{align}\label{ue28}
\frac{\partial^{r}}{\partial t^{r}}
e ^{-uL^t_2} =& \frac{(-1)^{k-1} (k-1)!}{2\pi i u^{k-1}}
\int_{\delta\cup \Delta}  e ^{-u\lambda}
\frac{\partial^{r}}{\partial t^{r}}(\lambda-L^t_2)^{-k} d \lambda.
\end{align}
We have the similar equation for
$\frac{\partial^{r}}{\partial t^{r}}(L^t_2e ^{-uL^t_2})$. Set
\be\label{ue29}
I_{k,r} = \Big \{ (\bk,\br)=(k_i,r_i)| \sum_{i=0}^j k_i =k+j, \sum_{i=1}^j r_i =r,\, \,   k_i, r_i \in \bN^*\Big \}.
\ee
Then  there exist $a ^{\bk}_{\br} \in \bR$ such that
\begin{align}\label{ue30}
& A^{\bk}_{\br}  (\lambda,t) = (\lambda-L^t_2)^{-k_0}
\frac{\partial^{r_1}L^t_2}{\partial t^{r_1}}  (\lambda-L^t_2)^{-k_1}
\cdots\frac{\partial^{r_j}L^t_2}{\partial t^{r_j}}  (\lambda-L^t_2)^{-k_j},\\
& \frac{\partial^{r}}{\partial t^{r}}
(\lambda-L^t_2)^{-k}=
\sum_{(\bk,\br)\in I_{k,r} }
 a ^{\bk}_{\br}  A ^{\bk}_{\br}  (\lambda,t).\nonumber
\end{align}

We claim that $A ^{\bk}_{\br}(\lambda,t)$ is well defined and
for any $m\in \bN$, $k>2(m+r+1)$, $Q,Q'\in \mQ^m$,
there exist $C>0$, $N\in \bN$ such that for $\lambda\in \delta\cup  \Delta$,
\begin{align}\label{1ue30}
\|Q  A ^{\bk}_{\br}(\lambda,t)Q' s\|_{t,0}
\leq C (1+|\lambda|)^N \sum_{|\beta|\leq 2r} \|Z^\beta s\|_{t,0}.
\end{align}

In fact, by (\ref{1ue1}), $\frac{\partial^{r}}{\partial t^{r}}L^t_2$ is
combination of $\frac{\partial^{r_1}}{\partial t^{r_1}}(g^{ij}(tZ))$
$(\frac{\partial^{r_2}}{\partial t^{r_2}}\nabla_{t,e_i})$
$(\frac{\partial^{r_3}}{\partial t^{r_3}}\nabla_{t,e_j})$,
 $\frac{\partial^{r_1}}{\partial t^{r_1}}(d(tZ))$,
$\frac{\partial^{r_1}}{\partial t^{r_1}}(d_{i}(tZ))$
$(\frac{\partial^{r_2}}{\partial t^{r_2}}\nabla_{t,e_i})$.
Now $\frac{\partial^{r_1}}{\partial t^{r_1}}(d(tZ))$
(resp. $\frac{\partial^{r_1}}{\partial t^{r_1}}\nabla_{t,e_i}$)
 ($r_1\geq 1$),
 are functions of  the type as
 $d'(tZ)Z^\beta$, $|\beta|\leq r_1$ (resp. $r_1+1$)
and $d'(Z)$ and its derivatives on $Z$ are  bounded smooth functions on $Z$.

Let $\cR'_t$ be the family of operators of the type
$$\cR'_{t} = \{ [f_{j_1} Q_{j_1},
[f_{j_2} Q_{j_2},\ldots [f_{j_l} Q_{j_l}, L^t_2]]\ldots ] \}$$
with $f_{j_i}$ smooth bounded (with its derivatives) functions
and $Q_{j_i}\in \{\nabla_{t,e_l}\}_{l=1}^{2n}$.

Now for the operator $A ^{\bk}_{\br}(\lambda,t)Q'$,
we will move first all the term $Z^\beta$ in $d'(tZ)Z^\beta$
as above to the right hand side of this operator, to do so,
we always use the commutator trick, i.e., each time,
we consider only the commutation for $Z_i$,
 not for $Z^\beta$ with $|\beta|>1$.
Then $A ^{\bk}_{\br}(\lambda,t)Q'$ is as the form
$\sum_{|\beta|\leq 2r} L^t_\beta Q''_\beta  Z^\beta$, and $Q''_\beta$ is obtained from
 $Q'$ and its commutation with $Z^\beta$.
Now we move all
the terms $\nabla_{t,e_i}$ in
$\frac{\partial^{r_j} L^t_2}{\partial t^{r_j}}$
to the right hand side of the operator $L^t_\beta$.
Then as in the proof of Theorem \ref{tu6},
we get finally  that $Q  A ^{\bk}_{\br}(\lambda,t)Q'$
is as the form  $\sum_{|\beta|\leq 2r} \cL^t_\beta Z^\beta$ where $\cL^t_\beta$
is a linear combination of operators of the form
\begin{align*}
Q (\lambda-L^t_2)^{-k'_0}R_1(\lambda-L^t_2)^{-k'_1}R_2
\cdots R_{l'}(\lambda-L^t_2)^{-k'_{l'}} Q''' Q'',
\end{align*}
with  $R_1, \cdots, R_{l'} \in \cR'_{t}$,  $Q'''\in \mQ^l$,
 $Q''\in   \mQ^m$, $|\beta|\leq 2 r$,
 and $Q''$ is obtained from $Q'$ and its commutation with $Z^\beta$.
By the argument as in \eqref{ue18} and  \eqref{ue20}, as $k>2(m+r+1)$,
we can split the above operator to two parts
\begin{align*}
&Q (\lambda-L^t_2)^{-k'_0}R_1(\lambda-L^t_2)^{-k'_1}R_2
\cdots R_{i}(\lambda-L^t_2)^{-k''_{i}};\\
&(\lambda-L^t_2)^{-(k'_{i}-k''_{i}) }\cdots
R_{l'}(\lambda-L^t_2)^{-k'_{l'}} Q''' Q'',
\end{align*}
and the $\|\quad \|^{0,0}_t$-norm of each part is bounded
by $C(1+|\lambda|^2)^N$. Thus
the proof of  (\ref{1ue30}) is complete.

\comment{
We claim that $A ^{\bk}_{\br}(\lambda,t)$ is well defined and
for any $l\in \bN$, there exist $C>0$, $N\in\bN$ such that for $\lambda \in \delta\cup \Delta$, $s\in
\cC^\infty_0(X_0,\bE_{x_0})$,
\begin{align}\label{1ue30}
\|A ^{\bk}_{\br}(\lambda,t) s\|_{t,l}
\leq C (1+|\lambda|)^N \sum_{|\alpha|\leq 2r} \|Z^\alpha s\|_{t,2r+l-k}.
\end{align}
In fact, by (\ref{1ue1}), $\frac{\partial^{r}}{\partial t^{r}}L^t_2$ is  combination of $\frac{\partial^{r_1}}{\partial
t^{r_1}}(g^{TX_0}_{ij}(tZ))
(\frac{\partial^{r_2}}{\partial t^{r_2}}\nabla_{t,e_i})
(\frac{\partial^{r_3}}{\partial t^{r_3}}\nabla_{t,e_j})$,
 $\frac{\partial^{r_1}}{\partial t^{r_1}}(b(tZ))$,
$\frac{\partial^{r_1}}{\partial t^{r_1}}(a_{i}(tZ))(\frac{\partial^{r_2}}{\partial t^{r_2}}\nabla_{t,e_i})$.
Now $\frac{\partial^{r_1}}{\partial t^{r_1}}(b(tZ))$
(resp. $\frac{\partial^{r_1}}{\partial t^{r_1}}\nabla_{t,e_i}$) ($r_1\geq 1$)
 are functions of  the type as
 $b'(tZ)Z^\beta$, $|\beta|\leq r_1$ (resp. $r_1+1$)
and $b'(Z)$ and its derivatives on $Z$ are  bounded smooth functions.
Thus by (\ref{ue12}), we get (\ref{1ue30}).
}

By (\ref{ue28}), (\ref{ue30}) and (\ref{1ue30}), we get the similar
estimates (\ref{ue15}) with $m'=C''=0$, (\ref{1ue27})
for $\frac{\partial^{r}}{\partial t^{r}}e ^{-uL^t_2}$,
$\frac{\partial^{r}}{\partial t^{r}}(L^t_2e ^{-uL^t_2})$.
Thus  we get (\ref{ue15}) for $m'=0$.

Finally, for $U$ a vector on $X$,
\begin{align}\label{0ue30}
& \nabla ^{\pi ^* \End (\bE)}_U e ^{-uL^t_2}
= \frac{(-1)^{k-1} (k-1)!}{2\pi i u^{k-1}}
\int_{\delta\cup \Delta}  e ^{-u\lambda}
\nabla ^{\pi ^* \End (\bE)}_U (\lambda-L^t_2)^{-k} d \lambda .
\end{align}
Now, by using the similar formula (\ref{ue30})  for
$\nabla ^{\pi ^* \End (\bE)} _U(\lambda-L^t_2)^{-k}$ by replacing
$\frac{\partial^{r_1}L^t_2}{\partial t^{r_1}}$ by
$\nabla ^{\pi ^* \End (\bE)}_U L^t_2$, and remark that
$\nabla ^{\pi ^* \End (\bE)}_U L^t_2$ is a differential operator
on $T_{x_0} X$ with the same structure  as $L^t_2$.
Then by the above argument, we get (\ref{ue15}) for $m'\geq 1$.
\end{proof}

Let $P_{0,t}$ be the orthogonal projection from
$\cC^\infty (X_0, \bE_{x_0})$ to the kernel of $L^t_2$ with respect to
$\left \langle \, ,\,\right \rangle_{t,0}$. Set
\begin{align}\label{1ue31}
F_u(L^t_2) = \frac{1}{2\pi i} \int_{\Delta}  e ^{-u\lambda}(\lambda-L^t_2)^{-1} d \lambda .
\end{align}
Let $P_{0,t}(Z,Z')$, $F_u(L^t_2)(Z,Z')$ be the smooth kernels of
 $P_{0,t},F_u(L^t_2)$ with respect to $dv_{TX}(Z')$. Then by (\ref{1u3}),
\begin{align}\label{1ue32}
F_u(L^t_2)=e ^{-uL^t_2}-P_{0,t} = \int_u^{+\infty} L^t_2e ^{-u_1L^t_2}du_1.
\end{align}
\begin{cor}\label{0tue8} With the notation in Theorem \ref{tue8},
\begin{multline}\label{1ue33}
\sup_{|\alpha|,|\alpha'|\leq m}
\Big |\frac{\partial^{|\alpha|+|\alpha'|}}
{\partial Z^{\alpha} {\partial Z'}^{\alpha'}}
\frac{\partial^{r}}{\partial t^{r}} F_u(L^t_2)
\left (Z, Z'\right )\Big |_{\cC^{m'}(X)}  \\
   \leq C (1+|Z|+|Z'|)^N
\exp ( -\frac{1}{8}\mu_0 u- \sqrt{C''\mu_0} |Z-Z'|).
\end{multline}
\end{cor}
\begin{proof} Note that $\frac{1}{8}\mu_0 u+\frac{2C''}{u} |Z-Z'|^2\geq \sqrt{C''\mu_0} |Z-Z'|$, thus
\begin{multline}\label{1ue34}
\int_u^{+\infty} e ^{ -\frac{1}{4}\mu_0 u_1- \frac{2C''}{u_1} |Z-Z'|^2}du_1
\leq e ^{-\sqrt{C''\mu_0}|Z-Z'|}\int_u^{+\infty} e ^{-\frac{1}{8}\mu_0 u_1}du_1 \\
 = \frac{8}{\mu_0}e ^{ -\frac{1}{8}\mu_0 u- \sqrt{C''\mu_0} |Z-Z'|}.
\end{multline}
By  (\ref{ue15}), (\ref{1ue32}), and (\ref{1ue34}), we get (\ref{1ue33}).
\end{proof}

\begin{rem}\label{0tue9} Under the condition of  Lindholm \cite{Li},
the metric on the trivial holomorphic  line bundle on $\bC^n$ is
$\|1\|=e ^{-\varphi/2}$. Now we use the unit section
$S_L= e ^{\varphi/2}1$ to trivialize this line bundle.
Then if $\varphi$ is $\cC^\infty$ and
$\frac{\partial ^\alpha}{\partial Z^\alpha} \varphi$ is bounded
for $|\alpha|\geq 3$,
 from (\ref{ue15}), (\ref{1ue32}), (\ref{1ue33}) with  $r=0$,
we can derive the off-diagonal estimate of the Bergman kernel  on
$\bC^n$. Actually, the $\cC^0$-estimate was obtained by  Lindholm
\cite[Prop. 9]{Li}.
\end{rem}

For $k$ large enough, set
\begin{align}\label{ue31}
& F_{r,u}
= \frac{(-1)^{k-1} (k-1)!}{2\pi i \, r! \, u^{k-1}}\int_{\Delta}
e ^{-u\lambda}   \sum_{(\bk,\br)\in I_{k,r} }
 a ^{\bk}_{\br}  A ^{\bk}_{\br}  (\lambda,0)d \lambda ,\\
&J_{r,u}= \frac{(-1)^{k-1} (k-1)!}{2\pi i \, r!\,  u^{k-1}}
\int_{\delta\cup \Delta}  e ^{-u\lambda}   \sum_{(\bk,\br)\in I_{k,r} }
 a ^{\bk}_{\br}  A ^{\bk}_{\br}  (\lambda,0)d \lambda ,\nonumber\\
&F_{r,u,t} = \frac{1}{r!}\frac{\partial^{r}}{\partial t^{r}}
F_u (L^t_2)- F_{r,u},
\quad J_{r,u,t} = \frac{1}{r!}\frac{\partial^{r}}{\partial t^{r}}
e ^{-uL^t_2}- J_{r,u}. \nonumber
\end{align}
Certainly, as $t\to 0$, the limit of  $\|\quad\|_{t,m}$ exists,
and we denote it by  $\|\quad\|_{0,m}$.
\begin{thm} \label{tue9} For any $r,k>0$, there exist $C>0$, $N\in \bN$ such that for
$t \in [0,1], \lambda \in \delta\cup \Delta$,
\begin{align}\label{ue32}
&  \left \|\Big (\frac{\partial^{r}L^t_2}{\partial t^{r}} -
\frac{\partial^{r}L^t_2}{\partial t^{r}}|_{t=0} \Big ) s \right \|_{t,-1}
\leq Ct \sum_{|\alpha|\leq r+3} \|Z^\alpha s\|_{0,1},\\
& \Big \|\Big (\frac{\partial^{r}}{\partial t^{r}} (\lambda-L^t_2)^{-k}
-\sum_{(\bk,\br)\in I_{k,r} }
 a ^{\bk}_{\br}  A ^{\bk}_{\br} (\lambda,0)\Big )s\Big \|_{0,0}
\leq C t  (1+|\lambda|^2)^N
\sum_{|\alpha|\leq 4r+3} \|Z^\alpha s\|_{0,0}.\nonumber
\end{align}
\end{thm}
\begin{proof}  Note that by (\ref{c37}), (\ref{u0}),
for $t\in [0,1]$, $k\geq 1$,
\begin{align}\label{1ue35}
\|s\|_{t,0}\leq C\|s \|_{0,0},\quad \|s\|_{t,k}\leq C \sum_{|\alpha|\leq k} \|Z^\alpha s\|_{0,k}.
\end{align}
 An application of Taylor expansion for
(\ref{1ue1}) leads to the following equation,
if $s,s'$ have compact support,
\be\label{ue33}
\Big | \left \langle  \Big (\frac{\partial^{r}L^t_2}{\partial t^{r}} -
\frac{\partial^{r}L^t_2}{\partial t^{r}} |_{t=0} \Big )s,s'\right \rangle_{0,0}\Big |
\leq C t \|s'\|_{t,1}\sum_{|\alpha|\leq r+3} \|Z^\alpha s\|_{0,1}.
\ee
Thus we get the first inequality of (\ref{ue32}).  Note that
\begin{align}\label{ue34}
& (\lambda-L^t_2)^{-1}- (\lambda-L^0_2)^{-1}
=(\lambda-L^t_2)^{-1}(L^t_2-L^0_2) (\lambda-L^0_2)^{-1}.
\end{align}
After taking the limit, we know that Theorems \ref{tu4}-\ref{tu6} still hold
for $t=0$.
From (\ref{ue2}), (\ref{ue33}) and (\ref{ue34}),
\begin{align}\label{0ue34}
& \left \|\left ((\lambda-L^t_2)^{-1}- (\lambda-L^0_2)^{-1} \right)s
 \right \|_{0,0}
\leq Ct  (1+|\lambda |^4)   \sum_{|\alpha|\leq 3} \|Z^\alpha s\|_{0,1}.
\end{align}
Now from the first inequality  of (\ref{ue32}) for $r=0$, (\ref{ue30})
and (\ref{0ue34}),  we get (\ref{ue32}).
\end{proof}

\begin{thm} \label{tue12}
There  exist $C>0$, $N\in \bN$
such that for $t \in ]0,1]$, $u\geq u_0$,
$q\in \bN$, $Z,Z'\in T_{x_0}X$, $|Z|, |Z'|\leq q$,
\begin{align}\label{ue42}
\Big |F_{r,u,t}(Z,Z')\Big |
\leq & C t^{1/2(2n+1)} (1+q)^N   e ^{-\frac{1}{8}\mu_0 u},\\
\Big |J_{r,u,t}(Z,Z')\Big |
\leq & C t^{1/2(2n+1)} (1+q)^N e ^{\frac{1}{2}\mu_0 u}.\nonumber
\end{align}
\end{thm}
\begin{proof}
Let $J^0_{x_0, q}$ be the vector space of square integrable sections of
 $\bE_{x_0}$ over $\{ Z\in T_{x_0}X, |Z|\leq q+1\}$. If $s\in J^0_{x_0, q}$,
 put $\|s\|^2_{(q)} = \int_{|Z|\leq q+1} |s|^2_{\bE_{x_0}} dv_{TX}(Z)$.
Let $\|A\|_{(q)}$ be the  operator norm of $A\in \cL(J^0_{x_0, q})$
with respect to $\|\quad\|_{(q)}$.  By (\ref{ue28}),
 (\ref{ue31}) and (\ref{ue32}), we get:
 There exist  $C>0, N\in \bN$
such that for $t \in ]0,1],  u\geq u_0$,
\begin{align}\label{ue39}
& \|F_{r,u,t}\|_{(q)}
\leq C t (1+q)^N   e ^{-\frac{1}{2}\mu_0 u},\\
& \|J_{r,u,t}\|_{(q)}
\leq C t (1+q)^N \ e ^{\frac{1}{4}\mu_0 u}.\nonumber
\end{align}

 Let $\phi : \bR\to [0,1]$ be a smooth function with compact support, equal $1$ near $0$, such that  $\int_{T_{x_0}X} \phi
(Z) dv_{TX}(Z)=1$. Take $\nu \in ]0,1]$. By the proof of Theorem \ref{tue8}, $F_{r,u}$ verifies the similar inequality as in
(\ref{1ue33}). Thus by (\ref{1ue33}), there exists $C>0$ such that if
$|Z|,|Z'|\leq q$, $U,U'\in \bE_{x_0}$,
\begin{multline}\label{ue43}
\Big |  \left \langle
F_{r,u,t} (Z,Z') U,U'     \right \rangle
-\int_{T_{x_0}X \times T_{x_0}X}  \left \langle
F_{r,u,t}(Z-W,Z'-W') U,U'     \right \rangle\\
  \frac{1}{\nu ^{4n}} \phi (W/\nu) \phi (W'/\nu) dv_{TX}(W)dv_{TX}(W')\Big |
\leq C \nu (1+q)^N e ^{ -\frac{1}{8}\mu_0 u} |U||U'|.
\end{multline}
On the other hand, by (\ref{ue39}),
\begin{multline}\label{ue44}
\Big |\int_{T_{x_0}X \times T_{x_0}X}  \left \langle
F_{r,u,t}  (Z-W,Z'-W') U,U'     \right \rangle\\
  \frac{1}{\nu ^{4n}} \phi (W/\nu) \phi (W'/\nu) dv_{TX}(W)dv_{TX}(W')\Big |
\leq C       t\frac{1}{\nu ^{2n}}  (1+q)^N e ^{-\frac{1}{2}\mu_0 u} |U||U'|.
\end{multline}
By taking $\nu = t^{1/2(2n+1)}$, we get (\ref{ue42}).
In the same way, we get  (\ref{ue42}) for $J_{r,u,t}$.
\end{proof}

\begin{thm} \label{tue14}  There exists $C''>0$ such that
for any $k,m,m'\in \bN$, there exist $N\in \bN$,
$C>0$  such that if
$t\in ]0,1], u\geq u_0$, $Z,Z'\in T_{x_0}X$,
 \begin{align}\label{0ue45}
&\sup_{|\alpha|,|\alpha'|\leq m}
\Big |\frac{\partial^{|\alpha|+|\alpha'|}}
{\partial Z^{\alpha} {\partial Z'}^{\alpha'}}\Big (F_u (L^t_2)
- \sum_{r=0}^k F_{r,u}t^r\Big ) (Z,Z')
\Big |_{\cC^{m'}(X)} \\
 &\hspace*{20mm}  \leq C  t^{k+1} (1+|Z|+|Z'|)^N
\exp ( -\frac{1}{8} \mu_0  u- \sqrt{C''\mu_0} |Z-Z'|),\nonumber\\
&\sup_{|\alpha|,|\alpha'|\leq m}
\Big |\frac{\partial^{|\alpha|+|\alpha'|}}
{\partial Z^{\alpha} {\partial Z'}^{\alpha'}}\Big (e ^{-uL^t_2} - \sum_{r=0}^k J_{r,u} t^r\Big ) (Z,Z')
\Big |_{\cC^{m'}(X)}\nonumber \\
&\hspace*{20mm}    \leq C  t^{k+1}  (1+|Z|+|Z'|)^N
\exp (\frac{1}{2} \mu_0 u- \frac{2C''}{u} |Z-Z'|^2).\nonumber
\end{align}
\end{thm}
\begin{proof} By (\ref{ue31}), and (\ref{ue42}),
 \begin{align}\label{0ue47}
&\frac{1}{r!}  \frac{\partial^{r}}{\partial t^{r}}
F_u (L^t_2) |_{t=0} = F_{r,u},\quad
 \frac{1}{r!}\frac{\partial^{r}}{\partial t^{r}} e ^{-uL^t_2}|_{t=0}
= J_{r,u}.
\end{align}
Now by Theorem \ref{tue8} and (\ref{ue31}), $ J_{r,u}, F_{r,u}$
have the same estimates as $\frac{\partial^{r}}{\partial t^{r}} e ^{-uL^t_2}$,
$\frac{\partial^{r}}{\partial t^{r}}F_u (L^t_2)$,
in (\ref{ue15}), (\ref{1ue33}).
Again from  (\ref{ue15}), (\ref{1ue33}), (\ref{ue31}),
and the Taylor expansion
$G(t)- \sum_{r=0}^k \frac{1}{r !} \frac{\partial ^r G}{\partial t^r}(0) t^r$
$= \frac{1}{k!}\int_0^t (t-t_0)^k
\frac{\partial ^{k+1} G}{\partial t^{k+1} }(t_0) dt_0$, we get (\ref{0ue45}).
\end{proof}

\subsection{Evaluation of $J_{r,u}$. } \label{s3.5}

For $u >0$, we will write $u\Delta_j$ for the rescaled simplex
$\{ (u_1, \cdots, u_j)|$ $0\leq u_1\leq u_2\leq \cdots\leq u_j\leq u\}$.
 By (\ref{c30}),
\begin{align}\label{ue48}
{\bf D}_t^2 = \mO_0^2
+\sum_{r=1}^\infty \sum_{r_1+r_2 =r} \mO_{r_1}\mO_{r_2}
t^r = L^0_2 +\sum_{r=1}^\infty \mQ _r t^r.
\end{align}
Set $\mathcal{J}= -2\pi \sqrt{-1} {\bf J}$. By  (\ref{0.1}), $\mathcal{J}\in \End (T^{(1,0)}X)$
is positive, and the $\mathcal{J}$ action on $TX$ is skew-symmetric.
We denote by $\det_{\bC}$ for the determinant function on the complex
 bundle $T^{(1,0)}X$.
We denote by
$|\mathcal{J}_{x_0}|=(\mathcal{J}^2_{x_0})^{1/2}$, and by $L^0_{2,\bC}$
  the restriction of $L^0_2$ on $\cC^\infty (\bR^{2n},\bC)$, then
by (\ref{b4}), (\ref{1ue1}),
\begin{align}\label{ue54}
&L^0_{2,\bC}=- \sum_j 
\Big(\nabla _{e_j} +\frac{1}{2}R^L_{x_0}(Z, e_j)\Big)^2   -\tau_{x_0},\\
&L^0_2 =L^0_{2,\bC} -2 \omega_{d,x_0}. \nonumber
\end{align}
 Let $e ^{-uL^0_{2,\bC}} (Z,Z')$,
$e ^{-uL^0_2} (Z,Z')$ be the smooth kernels of $e ^{-uL^0_{2,\bC}}$,
$e ^{-uL^0_2}$ with respect to $dv_{TX}(Z')$. Now from (\ref{ue54})
(cf. \cite[(6.37), (6.38)]{B90}),
\begin{align}\label{ue56}
e ^{-uL^0_{2,\bC}} (Z,Z')
=&\frac{1}{(2\pi)^{n}} {\det}_{\bC} \Big (\frac{\mathcal{J}_{x_0}}{1-e^{-2u\mathcal{J}_{x_0}}}\Big )
\exp  \Big (-\frac{1}{2} \left \langle
\frac{\mathcal{J}_{x_0}/2}{\tanh (u\mathcal{J}_{x_0})} Z,Z \right \rangle\\
&-\frac{1}{2} \left \langle \frac{\mathcal{J}_{x_0}/2}
{\tanh (u\mathcal{J}_{x_0})} Z',Z' \right \rangle
+ \left \langle \frac{\mathcal{J}_{x_0}/2}{\sinh (u\mathcal{J}_{x_0})}
e ^{u\mathcal{J}_{x_0}} Z,Z' \right \rangle
\Big ),\nonumber\\
e ^{-uL^0_2} (Z,Z')=& e ^{-uL^0_{2,\bC}} (Z,Z')e ^{2u \omega_{d,x_0}}.
\nonumber
\end{align}

\begin{thm}  \label{tue15}  For $r\geq 0$, we have
\begin{multline}\label{ue49}
J_{r,u}=\sum_{\sum_{i=1}^j r_i =r,\, r_i\geq 1} (-1)^j
\int_{u\Delta_j} e ^ {-(u-u_j)L^0_2} \mQ _{r_j}
e ^ {-(u_j-u_{j-1})L^0_2}\\
\cdots   \mQ _{r_1} e ^ {-u_1L^0_2} du_1\cdots du_j,
\end{multline}
where the product in the integrand is the convolution product. Moreover,
there exist  $J_{r, \beta, \beta'}(u)
\in \End(\Lambda (T^{*(0,1)}X)\otimes E)_{x_0}$
 smooth on $u\in ]0,+\infty [$ such that
\be\label{ue52}
J_{r,u}(Z,Z') = \sum_{|\beta|+|\beta'|\leq 3 r}
J_{r, \beta, \beta'}(u) Z^\beta {Z'}^{\beta'} e ^ {-u L^0_{2,\bC}}(Z,Z'),
\ee
and $\sum_{|\beta|+|\beta'|\leq 3 r} 
J_{r, \beta, \beta'}(u) Z^\beta {Z'}^{\beta'}$ as polynomial of $Z,Z'$
is even or odd according to whether $r$ is even or odd.
\end{thm}
\begin{proof}
We introduce an even extra-variable $\sigma$ such that $\sigma ^{r+1}=0$.
Set $[\quad]^{[r]}$  the coefficient of $\sigma ^r$,
 $L_\sigma= L^0_2 + \sum_{j=1}^r  \mQ _j \sigma ^j$.
From (\ref{ue31}), (\ref{0ue47}), we  know
\begin{align}\label{0ue51}
 J_{r,u}(Z,Z')= \frac{1}{r!} \frac{\partial ^r }{\partial t ^r} e ^{-uL^t_2}(Z,Z') |_{t=0}
=  [ e ^{-uL_\sigma}]^{[r]}(Z,Z').
\end{align}
Now from (\ref{0ue51}) and the Volterra expansion of
$e ^{-uL_\sigma}$ (cf. \cite[\S 2.4]{BeGeV}),
we get (\ref{ue49}).

We  prove (\ref{ue52}) by iteration.
By (\ref{ue56}), (\ref{ue49}) and Theorem \ref{t3.3},
 we immediately derive (\ref{ue52}).
By the iteration, (\ref{ue56}) and  Theorem \ref{t3.3},
the polynomial of $Z,Z'$ has the same parity with $r$.
\end{proof}

\subsection{Proof of Theorems \ref{t0.1}, \ref{t0.2}}
\label{s3.6}

By (\ref{1ue32}), (\ref{0ue45}),  for any $u>0$ fixed,
 there exists $C_u>0$ such that for
$t=\frac{1}{\sqrt{p}}$, $Z,Z'\in T_{x_0}X$, $x_0\in X$, we have
\begin{multline}\label{1ue54}
\sup_{|\alpha|,|\alpha'|\leq m}\left |\frac{\partial^{|\alpha|+|\alpha'|}}
{\partial Z^{\alpha} {\partial Z'}^{\alpha'}} \Big (P_{0,t}
- \sum_{r=0}^k t^r(J_{r,u}- F_{r,u})\Big)(Z,Z')\right |_{\cC^{m'}(X)}\\
\leq C_u  t^{k+1}  (1+|Z|+|Z'|)^N
\exp (- \sqrt{C''\mu_0 } |Z-Z'|),
\end{multline}
Set
\begin{align}\label{1ue53}
P^{(r)}=J_{r,u}- F_{r,u}.
\end{align}
Then  $P^{(r)}$ does not depend on $u>0$ by (\ref{1ue54}), as $P_{0,t}$
does not depend on $u$.
Moreover, by taken the limit of (\ref{1ue33}) as $t\to 0$,
\begin{align}\label{1ue56}
\Big | F_{r,u}(Z,Z')\Big |_{\cC^m(X)}
\leq C (1+|Z|+|Z'|)^N \exp ( -\frac{1}{8} \mu_0  u- \sqrt{C''\mu_0} |Z-Z'|).
\end{align}
Thus
\begin{align}\label{1ue57}
J_{r,u}(Z,Z') = P^{(r)}(Z,Z')+F_{r,u}(Z,Z')= P^{(r)}(Z,Z')+\cO( e^ {-\frac{1}{8} \mu_0  u}),
\end{align}
uniformly on any compact set of $T_{x_0}X\times T_{x_0}X$.

%We denote by
%$|\mathcal{J}_{x_0}|=(\mathcal{J}^2_{x_0})^{1/2}$, and by $L^0_{2,\bC}$
%  the restriction of $L^0_2$ on $\cC^\infty (\bR^{2n},\bC)$, then
%\begin{align*}
%L^0_{2,\bC}=- \sum_j
%\Big(\nabla _{e_j} +\frac{1}{2}R^L_{x_0}(Z, e_j)\Big)^2   -\tau_{x_0}.
%\end{align*}
Let $P(Z,Z')$ be the Bergman kernel of $L^0_{2,\bC}$ in \eqref{ue54},
 i.e. the smooth kernel of the orthogonal projection from
$L^2 (\bR^{2n},\bC)$ on $\Ker L^0_{2,\bC}$.
Then for $Z,Z'\in T_{x_0}X$,
\begin{align}\label{ue62}
P(Z,Z')  =\frac{\det_{\bC}{\mathcal{J}_{x_0}}}{(2\pi)^n}
\exp\Big (- \frac{1}{4} \left \langle |\mathcal{J}_{x_0}|(Z-Z'),(Z-Z') \right \rangle
+\frac{1}{2} \left \langle \mathcal{J}_{x_0} Z,Z' \right \rangle\Big ).
\end{align}
Now $e^{ u\mathcal{J}_{x_0}} = \cosh (u|\mathcal{J}_{x_0}|)+ \sinh
(u|\mathcal{J}_{x_0}|)\frac{\mathcal{J}_{x_0}}{|\mathcal{J}_{x_0}|}$, thus
$\frac{\mathcal{J}_{x_0}/2}{\sinh(u\mathcal{J}_{x_0})}e^{ u\mathcal{J}_{x_0}}= \frac{1}{2} (|\mathcal{J}_{x_0}|+
\mathcal{J}_{x_0}) +\cO(e ^{-2u |\mathcal{J}_{x_0}|})$.
From (\ref{ue56}), and (\ref{ue49}), we get as $u\to \infty$,
\begin{align}\label{ue63}
&J_{0,u}(Z, Z') =e ^{-uL^0_2} (Z, Z')
=  P(Z,Z')I_{\bC\otimes E} + \cO( e^ {-\mu_0  u}),\\
& P^{(0)}(Z,Z')=P(Z,Z')I_{\bC\otimes E}.\nonumber
\end{align}
uniformly on any compact set of $T_{x_0}X\times T_{x_0}X$.
From (\ref{ue52}), (\ref{1ue57}), and (\ref{ue63}),
we know that as $u\to \infty$,
\begin{align}\label{1ue58}
J_{r,\beta,\beta'}(u)=J_{r,\beta,\beta'}(\infty)+ \cO( e^ {-\frac{1}{8} \mu_0  u}).
\end{align}
and by (\ref{1ue57}), (\ref{ue63}) and (\ref{1ue58}),
\begin{align}\label{1ue59}
& P^{(r)}(Z,Z')=J_{r,\infty}(Z,Z')
= \sum_{\beta,\beta'}J_{r,\beta,\beta'}(\infty)
Z^\beta {Z'}^{\beta'}P(Z,Z').
\end{align}

Note that in (\ref{c22}), $\kappa (Z) =(\det g_{ij}(Z))^{1/2}
= (\det (\theta_i^k \theta_j^k))^{1/2}$.
By (\ref{c27}), for $Z,Z'\in T_{x_0}X$,
\begin{align}\label{1ue60}
&P^0_p(Z,Z')= p^n P_{0,t}(Z/t,Z'/t) \kappa ^{-1}(Z'), \\
&\exp(-\frac{u}{p} D^{X_0,2}_p)(Z,Z')
= p^n  e ^{-uL^t_2} (Z/t,Z'/t)\kappa ^{-1}(Z').\nonumber
\end{align}

We now observe that, as a consequence of  (\ref{1ue54}) and
(\ref{1ue60}), we obtain the following important estimate.
\begin{thm} \label{tue17}
For any $k,m,m'\in \bN$, there exist $N\in \bN, C>0$ such that for
$\alpha, \alpha' \in \bN^{n}$, $|\alpha|+|\alpha'|\leq m$,
$Z,Z'\in T_{x_0}X$, $|Z|, |Z'|\leq  \var$, $x_0\in X$, $p\geq 1$,
\begin{multline}\label{ue66}
\left |\frac{\partial^{|\alpha|+|\alpha'|}}
{\partial Z^{\alpha} {\partial Z'}^{\alpha'}}
\left (\frac{1}{p^n}  P_p^0(Z,Z')
-\sum_{r=0}^k  P^{(r)} (\sqrt{p} Z,\sqrt{p} Z')\kappa ^{-1}(Z')
p^{-r/2}\right )\right |_{\cC^{m'}(X)}\\
\leq C  p^{-(k+1-m)/2}  (1+|\sqrt{p} Z|+|\sqrt{p} Z'|)^N
\exp (- \sqrt{C''\mu_0 } \sqrt{p} |Z-Z'|).
\end{multline}
\end{thm}

By (\ref{1c19}) and Theorem \ref{tue17},
we obtain the following full off-diagonal expansion
for the Bergman kernel on $X$.

{\bf Theorem 3.18$^\prime$.} \label{atue17}
{\em With the notation in Theorem \ref{tue17},}
\begin{multline}\label{aue66}
\left |\frac{\partial^{|\alpha|+|\alpha'|}}
{\partial Z^{\alpha} {\partial Z'}^{\alpha'}}
\left (\frac{1}{p^n}  P_p(Z,Z')
-\sum_{r=0}^k  P^{(r)} (\sqrt{p} Z,\sqrt{p} Z')\kappa ^{-1}(Z')
p^{-r/2}\right )\right |_{\cC^{m'}(X)}\\
\leq C  p^{-(k+1-m)/2}  (1+|\sqrt{p} Z|+|\sqrt{p} Z'|)^N
\exp (- \sqrt{C''\mu_0 } \sqrt{p} |Z-Z'|)+ \cO(p^{-\infty}).
\end{multline}
The term $\cO(p^{-\infty})$ means that for any $l,l_1\in \bN$,
there exists $C_{l,l_1}>0$ such that its $\cC^{l_1}$-norm is dominated
by $C_{l,l_1} p^{-l}$.

From Theorem \ref{tue15}, we know that $J_{r,u}(0,0)=0$ for $r$ odd. Thus from (\ref{1ue57}), $P^{(r)}(0,0)=0$ for $r$
odd. Thus from
 (\ref{ue66}), for $Z=Z'=0$, $m=0$, we get
\begin{align}\label{ue67}
\Big  |\frac{1}{p^n}  P_p(x_0,x_0 ) -\sum_{r=0}^k
 P^{(2r)}(0,0)p^{-r}\Big  |_{\cC^{m'}(X)}
\leq C p^{-k-1}.
\end{align}
From  (\ref{ue63}),
\begin{align}\label{ue68}
P^{(0)}(0,0)=P(0,0)I_{\bC\otimes E}=
 (\det {\bf J})^{1/2} I_{\bC\otimes E}.
\end{align}
Moreover, from Theorems \ref{t3.3}, \ref{tue15}, (\ref{ue48}), we
deduce the desired property on $b_{r}$ in  Theorem \ref{t0.1}. To
get the last part of Theorem \ref{t0.1}, we notice that the
constants in Theorems \ref{tue8} and \ref{tue12} will be uniform
bounded under our condition, thus we can take $C_{k,\, l}$ in
(\ref{0.6}) independent of $g^{TX}$. Thus we have proved Theorem
\ref{t0.1}.

From Proposition \ref{t3.2},  we know that for
any $u>0$ fixed, for any $l\in \bN$,
 there exists $C>0$ such that for $Z,Z'\in T_{x_0}X$,
  $|Z|, |Z'|\leq  \var$, $x_0\in X$,
\begin{align}\label{ue69}
\Big |\Big ( \exp(-\frac{u}{p} D^{2}_p)-  \exp(-\frac{u}{p} D^{X_0,2}_p)\Big )(Z,Z') \Big |_{\cC^{m'}(X)} \leq C p^{-l}.
\end{align}
Thus from Theorem \ref{tue15}, (\ref{0ue45}), (\ref{1ue60}), and (\ref{ue69}),
we get
\begin{align}\label{1ue69}
\Big |\frac{1}{p^n} \exp(-\frac{u}{p} D^{2}_p)(x_0,x_0)
-  \sum_{r=0}^k  J_{2r,u}(0,0)p^{-r}\Big |_{\cC^{m'}(X)}
 \leq C p^{-k-1}.
\end{align}
Hence we have (\ref{0.4}) and at $x_0$,
\begin{align}\label{ue70}
b_{r,u} = J_{2r,u}(0,0).
\end{align}
 Now, from (\ref{ue56}), (\ref{1ue57}), (\ref{ue67}), and (\ref{ue70}),
we deduce Theorem \ref{t0.2}.\\

From our proof of Theorems \ref{t0.1}, \ref{t0.2},
we also obtain a method to compute the coefficients.
Namely, we compute first the heat kernel expansion of
$\exp(-\frac{u}{p}D^2_p)(x,x)$ when $p\to \infty$ by
$\sum_{r=0}^j b_{r,u}(x) p^{n-r}$ (cf. (\ref{1ue69})),  then let $u\to \infty$,
we get the corresponding coefficients of the expansion
of $\frac{1}{p^n}  P_p(x,x)$.
As an example, we will calculate $b_1$ in the next section.

In practice, we choose $\{ w_i\}_{i=1}^n$
an orthonormal basis of $T^{(1,0)}_{x_0} X$, such that
\be\label{0ue52}
\mathcal{J}_{x_0}=
{\rm diag} (a_1 (x_0), \cdots, a_n(x_0))\in {\rm End} (T^{(1,0)}_{x_0} X),
\ee
with $0< a_1(x_0)\leq a_2(x_0)\leq \cdots \leq  a_n(x_0)$,
and let $\{w^j\}_{j=1}^n$ be its dual basis.
Then
$e_{2j-1}=\tfrac{1}{\sqrt{2}}(w_j+\overline{w}_j)$ and
$e_{2j}=\tfrac{\sqrt{-1}}{\sqrt{2}}(w_j-\overline{w}_j)\,,
 j=1,\dotsc,n\, $
forms an orthonormal basis of $T_{x_0}X$.
In the coordinate induced by $\{ e_i\}$ as above,
all even function $g(\mathcal{J}_{x_0})$ of $\mathcal{J}_{x_0}$ is diagonal,
and $g(\mathcal{J}_{x_0})=g(|\mathcal{J}_{x_0}|)$.

\section{Applications} \label{s4}

This section is organized as follows. In Section \ref{s4.1},
we calculate the coefficient $b_1$ in Theorem \ref{t0.1}
when the manifold is K\"ahler.
In Section \ref{s4.2}, we extend Theorem \ref{t0.1} to the orbifold case.
Again the finite propagation speed
 allows us to localize the problem
which was also used in \cite{M}.
\comment{
Contrast to the usual index theorem for orbifolds
\cite{K}, \cite{M}, we don't see the singularity term
in the asymptotic expansion, but the finite propagation speed
still allows us to localize the problem
which was also used in \cite{M}.
}

\subsection{K\"ahler case}\label{s4.1}

In this Section, we assume that $(X,\omega)$ is K\"ahler and ${\bf J}=J$,
and the vector bundles $E,L$ are holomorphic on $X$.
Then $a_j(x)=2\pi$ for $j\in \{1,\cdots, n\}$ in (\ref{0ue52}).
Note that for $\{w_j\}$ (resp. $\{e_j\}$) an orthonormal basis of $T^{(1,0)}X$
(resp $TX$), the scalar curvature $r^X$  of $(X, g^{TX})$ is given by
\be\label{d0}
r^X = -\sum_{jk}\left \langle R^{TX} (e_j,e_k)e_j,e_k \right \rangle
=2 \sum_{jk}\left \langle
R^{TX} (w_j,\overline{w}_j)w_k,\overline{w}_k\right \rangle.
\ee

Now the Levi-Civita connection $\nabla ^{TX}$ preserves $T^{(1,0)}X$
and $T^{(0,1)}X$, and
$\nabla ^{T^{(1,0)}X} = P^{T^{(1,0)}X}\nabla ^{TX} P^{T^{(1,0)}X}$
is the holomorphic Hermitian connection on $T^{(1,0)}X$.
In this situation,  the Clifford connection
$\nabla ^{{\rm Cliff}}$ on $\Lambda (T^{*(0,1)}X)$ is
$\nabla ^{ \Lambda (T^{*(0,1)}X)}$, the natural  connection induced by
$\nabla ^{T^{(1,0)}X}$.
Let $\overline{\partial} ^{L^p\otimes E,*}$  be the formal adjoint of
the Dolbeault operator
$\overline{\partial} ^{L^p\otimes E}$ on
$\Omega ^{0,\bullet}(X, L^p\otimes E)$.
Then  the operator $D_p$ in (\ref{defDirac})
is $D_p = \sqrt{2}( \overline{\partial} ^{L^p\otimes E}
+ \overline{\partial} ^{L^p\otimes E,*})$.
 Note that $D_p^2$ preserves the $\bZ$-grading of
  $\Omega ^{0,\bullet}(X, L^p\otimes E)$.
Let $D^2_{p,i}= D^2_p|_{\Omega ^{0,i}(X, L^p\otimes E)}$, then for
$p$ large enough,
\be \label{d1}
\Ker D_p = \Ker D^2_{p,0} = H^0 (X,L^p\otimes E).
\ee

By (\ref{d1}), $B_p(x)\in \End (E)$ and
 we only need to do the computation for $D^2_{p,0}$.
In what follows, we compute everything on $\cC^\infty (X, L^p\otimes E)$.
Especially, $\mQ_r$ in (\ref{ue48}) takes value in $\End (E)$.
Now, we replace $X$ by $\bR^{2n} \simeq T_{x_0}X$ as in Section \ref{s3.2},
and we use the notation therein.
We denote by $(g^{ij}(Z))$ the inverse of
the matrix $(g_{ij}(Z))=(g^{TX}_{ij}(Z))$.
 Let $\Delta^{TX} = \sum_i \frac{\partial ^2}{\partial Z_i^2}$
be the standard Euclidean Laplacian on $T_{x_0}X$ with
respect to the metric $g^{T_{x_0}X}$. Then by (\ref{0c25}), (\ref{0c31}),
\begin{align}\label{0d0}
g_{ij}(Z)=\sum_k \theta ^k_i \theta ^k_j (Z) =
\delta_{ij} +  \frac{1}{3}
\left \langle R^{TX}_{x_0} (\mR,e_i) \mR, e_j\right \rangle
 + \cO(|Z|^3).
\end{align}
\begin{thm} \label{t4.1}
\begin{align} \label{0d1}
&\mQ_0= -\Delta ^{TX} +\pi ^2 |Z|^2 - 2\pi n
+ 2 \sqrt{-1}\pi \nabla_{J\mR} ,
\quad  \mQ_1=0,\\
&\mQ_2= \sum_j \Big (\frac{2}{3}  \left \langle R^{TX}_{x_0} (\mR,e_i) e_i,
e_j\right \rangle -\frac{\sqrt{-1} \pi }{2 }
\left \langle R^{TX}_{x_0} (\mR,J\mR) \mR,e_j \right \rangle
 -R^E_{x_0} (\mR, e_j)\Big)\nabla_{e_j} \nonumber \\
&\hspace*{3mm}- \sqrt{-1} \sum_j \Big( \frac{1}{2}  R^E_{x_0} (e_j, Je_j)
+\frac{\pi}{2}
 \left \langle R^{TX}_{x_0} (\mR,e_j)\mR, Je_j\right \rangle \Big)
+ \pi \sqrt{-1}  R^E_{x_0} (\mR, J\mR)\nonumber \\
&\hspace*{3mm}- \frac{\pi ^2}{6}
\left \langle R^{TX}_{x_0} (\mR,J\mR) \mR,J\mR \right \rangle
+\frac{1}{3} \sum_{ij} \left \langle R^{TX}_{x_0} (\mR,e_i) \mR, e_j\right \rangle \nabla_{e_i} \nabla_{e_j} .\nonumber
\end{align}
\end{thm}
\begin{proof} Let $ \Gamma _{ij}^l$ be the connection form of $\nabla ^{TX}$ with respect to the basis $\{e_i\}$, then
$(\nabla ^{TX}_{e_i}e_j)(Z) = \Gamma _{ij}^l (Z) e_l$. By (\ref{0d0}),
\begin{align}\label{0d2}
\Gamma _{ij}^l (Z)& =  \frac{1}{2} \sum_k g^{lk} (\partial_i g_{jk}
+ \partial_j g_{ik}-\partial_k g_{ij})(Z)\\
&= \frac{1}{3}\Big  [ \left \langle R^{TX}_{x_0} (\mR, e_j) e_i, e_l\right \rangle _{x_0}
+ \left \langle R^{TX}_{x_0} (\mR, e_i) e_j, e_l\right \rangle_{x_0}\Big ]
 + \cO(|Z|^2). \nonumber
\end{align}
Observe that $J$ is parallel with respect to $\nabla ^{TX}$, thus $\left \langle J\wi{e}_i,\wi{e}_j \right \rangle_{Z}=\left
\langle Je_i,e_j \right \rangle_{x_0}$.
 By (\ref{0.1}), (\ref{0c25}), and (\ref{0c31}),
\begin{align}\label{0d3}
&\frac{\sqrt{-1}}{2 \pi} R^L(e_k,e_l)(Z)= \sum_{ij}\theta ^i_k(Z)\theta ^j_l(Z)
\left \langle J\wi{e}_i,\wi{e}_j \right \rangle_{Z}\\
&= \left \langle Je_k,e_l \right \rangle_{x_0}
-\frac{1}{6}\left \langle R^{TX}_{x_0} (\mR, e_k)\mR, Je_l\right \rangle_{x_0}
 +\frac{1}{6}\left \langle R^{TX}_{x_0} (\mR, Je_k)\mR, e_l\right \rangle_{x_0}
  + \cO(|Z|^3). \nonumber
\end{align}

By (\ref{c27}), (\ref{c35}) and (\ref{0d3}), for $t=\frac{1}{\sqrt{p}}$,
we get
\begin{multline}\label{d8}
\nabla_{t, e_i}|_{Z} =  tS_t^{-1}\nabla ^{L^p\otimes E}_{e_i}S_t|_{Z}
= \nabla_{e_i}
+ \frac{1}{t} \Gamma ^L (e_i)(tZ) + t \Gamma ^E (e_i) (tZ)\\
= \nabla_{e_i} -\sqrt{-1} \pi \left \langle J\mR,e_i\right \rangle
 -\frac{\sqrt{-1}\pi}{12}t^2 \left \langle R^{TX}_{x_0} (\mR, J\mR ) \mR, e_i\right \rangle + \frac{t^2}{2} R^E_{x_0}
(\mR,e_i)
+ \cO(t^3).
\end{multline}
By a direct calculation (\ref{1ue1}) or by Lichnerowicz formula in
\cite[Proposition 1.2]{B}, we know
\be\label{d4}\quad
D^2_{p,0} = - \sum_{ij}g^{ij} [\nabla ^{L^p\otimes E}_{e_i}\nabla ^{L^p\otimes E}_{e_j}
- \Gamma _{ij}^l \nabla ^{L^p\otimes E}_{e_l} ]
- \frac{\sqrt{-1}}{2}\sum_{i} R^E (\wi{e}_i, J\wi{e}_i)
-2 \pi n p.
\ee
Thus from  (\ref{c29}), (\ref{0d0}), (\ref{d8}), and (\ref{d4}),
\begin{align}\label{d9}
&{\bf D}_{t,0}^2= S_t^{-1}  t^2 D_{p,0}^2  S_t\\
&= -\sum_{ij}g^{ij} (tZ) \Big [ \nabla_{t, e_i}\nabla_{t, e_j}
- t \Gamma _{ij}^l(t\cdot) \nabla_{t, e_l} \Big ](Z)
-  \frac{\sqrt{-1}}{2} t^2 \sum_{i}R^E (\wi{e}_i, J\wi{e}_i) (tZ) -2 \pi n \nonumber \\
&= -\sum_{ij} \Big ( \delta_{ij} - \frac{t^2}{3}
\left \langle R^{TX}_{x_0} (\mR,e_i) \mR, e_j\right \rangle
 + \cO(t^3)\Big ) \nonumber \\
&  \Big \{ \Big (\nabla_{e_i} -\sqrt{-1} \pi \left \langle J\mR,e_i\right \rangle -\frac{\sqrt{-1}\pi}{12}t^2 \left \langle
R^{TX}_{x_0} (\mR, J\mR ) \mR, e_i\right \rangle
+ \frac{t^2}{2} R^E_{x_0} (\mR,e_i)
+ \cO(t^3)\Big ) \nonumber \\
& \Big (\nabla_{e_j} -\sqrt{-1} \pi \left \langle J\mR,e_j\right \rangle
-\frac{\sqrt{-1}\pi}{12}t^2
\left \langle R^{TX}_{x_0} (\mR, J\mR ) \mR, e_j\right \rangle
+ \frac{t^2}{2} R^E _{x_0}(\mR,e_j) + \cO(t^3)\Big )   \nonumber \\
&- t \Gamma ^l_{ij}(tZ) \Big (\nabla_{e_l}
-\sqrt{-1} \pi \left \langle J\mR,e_l\right \rangle
+ \cO(t^2)\Big ) \Big \}  \nonumber \\
&- \frac{\sqrt{-1}}{2} t^2 \sum_{i}R^E_{x_0} (e_i, Je_i) -2\pi n
+ \cO(t^3). \nonumber
\end{align}

From  (\ref{0d2}), (\ref{d9}) and the fact that $R^{TX}$ is a (1,1)-form,
we derive (\ref{0d1}).
\end{proof}

\comment{The following result was obtained in \cite{Catlin}, \cite{Lu}, \cite{Wang1},
\cite{Zelditch} in various generality.
\begin{thm} \label{t4.2}If $(X,\omega)$ is K\"ahler and ${\bf J}=J$,
then there exist smooth coefficients
$b_j(x)\in \End (E)_x$
with
\be \label{d10}
b_0=\Id _E, \quad b_1 = \frac{1}{4\pi}
\Big [ \sqrt{-1} \sum_{i} R^E (e_i, Je_i) +\frac{1}{2}r^X \Id_E \Big ].
\ee
such that for any $k,l\in \bN$, there exists
$C_{k,l}>0$ such that for any $x\in X$, $p\in \bN$, we have (\ref{0.3}).
Moreover, $b_r(x)$ are polynomials in $R^{TX}$,
$R^E$  and their derivatives with order
$\leq 2r-1$  at $x$.
\end{thm}}

\begin{proof}[Proof of Theorem \ref{t0.3}] From (\ref{ue56}), and (\ref{0d1}),
\begin{multline}\label{0d10}
e ^{-uL^0_2} (Z,Z') = \frac{1}{(1-e ^{-4\pi u})^n}
\exp \Big (-\frac{\pi (|Z|^2+|Z'|^2)}{2\tanh (2\pi u)}
+ \frac{\pi}{\sinh (2\pi u)} \left \langle
e ^{-2\sqrt{-1} \pi u J} Z,Z' \right \rangle \Big ).
\end{multline}
By (\ref{ue49}), (\ref{0d1}), (\ref{0d10}), $J_{1,u}(Z,Z')=0$, and
\begin{multline}\label{0d11}
J_{2,u}(0,0) = - \int^u_0 du_1 \int_{\bR ^{2n}}
\frac{1}{(1-e ^{-4\pi u_1})^n(1-e ^{-4\pi (u-u_1)})^n} \\
 \exp \Big (-\frac{\pi |Z|^2 }{2\tanh (2\pi (u-u_1))} \Big )\mQ_2(Z)
\exp \Big (-\frac{\pi |Z|^2 }{2\tanh (2\pi u_1)} \Big ).
\end{multline}
By (\ref{0d1}),
\begin{multline}\label{0d12}
\mQ_2(Z) \exp \Big (\frac{-\pi |Z|^2 }{2\tanh (2\pi u_1)} \Big )=
\Big \{\pi \sqrt{-1}  R^E_{x_0} (\mR, J\mR)
-\frac{\pi ^2}{6}\left \langle R^{TX}_{x_0} (\mR,J\mR) \mR,J\mR \right \rangle
 \\
-\frac{\sqrt{-1}}{2} \pi\sum_i \left \langle R^{TX}_{x_0} (\mR,e_i)\mR,
Je_i \right \rangle+ \frac{\pi}{3 \tanh (2\pi u_1)} \sum_i
\left \langle R^{TX}_{x_0} (\mR,e_i)\mR, e_i \right \rangle\\
- \frac{\sqrt{-1}}{2} \sum_i R^E_{x_0} (e_i, Je_i)\Big \}
\exp \Big (\frac{-\pi |Z|^2 }{2\tanh (2\pi u_1)} \Big ).
\end{multline}
Now $\int_{-\infty}^{+\infty} x^2 e ^{-x^2/2} dx
= \int_{-\infty}^{+\infty}  e ^{-x^2/2} dx = \sqrt{2\pi}$, and
 $\int_{-\infty}^{+\infty} x^4 e ^{-x^2/2} dx= 3\sqrt{2\pi}$. Thus
\begin{multline}\label{0d13}
\int_{\bR^{2n}} \left \langle R^{TX}_{x_0} (\mR,J\mR) \mR,J\mR \right \rangle
\exp(-\frac{|Z|^2}{2})
= (2\pi)^n \sum_{jk}
\Big [ \left \langle
R^{TX}_{x_0} (e_j,Je_j)e_k,Je_k \right \rangle\\
+\left \langle R^{TX}_{x_0} (e_j,Je_k)e_j,Je_k \right \rangle
+\left \langle R^{TX}_{x_0} (e_j,Je_k)e_k,Je_j \right \rangle\Big ]\\
= - (2\pi)^n \times4 r^X_{x_0} .
\end{multline}

Set $c(u_1)=\frac{\sinh (2\pi (u-u_1))\sinh (2\pi u_1)}{\sinh (2\pi u)}$.
 Then from (\ref{0d11})-(\ref{0d13}), we get
\begin{multline}\label{0d14}
J_{2,u}(0,0) = - \int^u_0 \frac{ du_1 }{(1-e ^{-4\pi u})^n}
\Big [ (c(u_1)-\frac{1}{2} )\sqrt{-1} \sum_i R^E_{x_0} (e_i, Je_i)\\
+ \frac{1}{3} \Big ( \frac{c(u_1)}{\tanh(2\pi u_1)}-2 c(u_1)^2 \Big) r^X_{x_0} \Big ]\\
= \frac{-1}{(1-e ^{-4\pi u})^n} \left\{
\Big [\Big(\frac{1}{\tanh (2\pi u)} -1\Big) \frac{u}{2} -\frac{1}{4\pi}\Big]\sqrt{-1} \sum_i R^E_{x_0} (e_i, Je_i)\right.\\
\left. + \frac{1}{3} \left[\frac{u}{2}- \frac{u}{2\tanh ^2 (2\pi u)}
- \frac{2}{\sinh^2(2\pi u)} \Big( \frac{-3}{32\pi } \sinh (4\pi u) +\frac{u}{8}\Big)\right]r^X_{x_0} \right \}.
\end{multline}
Thus by (\ref{0.5}), and (\ref{ue70}),
\be\label{0d15}
b_1 = \lim_{u\to \infty} J_{2,u}(0,0) = \frac{1}{4\pi}
\Big [ \sqrt{-1} \sum_i R^E (e_i, Je_i) +\frac{1}{2}r^X \Id_E \Big ].
\ee
From Theorem \ref{t0.1} and (\ref{0d15}),
the proof of Theorem \ref{t0.3} is completed.
\end{proof}

\subsection{Orbifold case}\label{s4.2}
Let $(X,\om)$ be a compact symplectic orbifold  of real dimension $2n$
with singular set $X'$.
By definition, for any $x\in X$, there  exists a small neighborhood
$U_x\subset X$,
 a finite group $G_x$ acting linearly on $\bR^{2n}$, and
$ \widetilde{U}_x\subset \bR^{2n}$ an $G_x$-open set such that
$\widetilde{U}_x \stackrel{\tau_x}{\to} \widetilde{U}_x/G_x=U_x$
and $\{0\}=\tau ^{-1}_x(x) \in \widetilde{U}_x$. We will use $\widetilde{z}$ to denote the point in $\widetilde{U}_x$
representing $z\in {U}_x$. Let $\Sigma X=\{(x,(h_x^j))|x\in X, G_x \neq 1,
(h_x^j)$ runs over the conjugacy  classes in $G_x\}$.
Then  $\Sigma X $ has a natural orbifold structure defined by (cf. \cite{K})
\begin{eqnarray}\label{e0}
\Big \{(Z_{G_x}(h_x^j)/K_x^j, \widetilde{U}_x^{h_x^j})
\rightarrow  \widetilde{U}_x^{h_x^j}/Z_{G_x}(h_x^j) \Big \}_{(x, U_x,j)}.
\end{eqnarray}
Here $\widetilde{U}_x^{h_x^j}$ is the  fixed point set of $h_x^j$
over $\widetilde{U}_x$, $Z_{G_x}(h_x^j)$ is the centralizer of $h_x^j$
in $G_x$, and $K_x^j$ is the kernel of the representation
$Z_{G_x}(h_x^j) \rightarrow$ Diffeo $(\widetilde{U}_x^{h_x^j})$.
The number $|K_x^j|$ is locally constant on $\Sigma X$
and we call it as the multiplicity $m_i$ of each connected component
$X_i$ of $X\cup \Sigma X$.

An orbifold vector bundle $E$ on  an orbifold $X$ means
that for any $x \in X$,  there exists
$ \widetilde{p}_{U_x}: \widetilde{E}_{U_x}
\rightarrow \widetilde{U_x}$
 a $G_{U_x}^{E}$-equivariant vector bundle and $(G_{U_x}^{E},
\widetilde{E}_{U_x})$ (resp. $(G_{U_x}^{E}/K_{U_x}, \widetilde{U}_x)$,
$K_{U_x} = {\rm Ker} (G_{U_x}^{E}\rightarrow {\rm Diffeo} (\widetilde{U}_x)))$
is the orbifold structure of $E$ (resp. $X$).
Set $\widetilde{E}_{U_x}^{\mathrm{pr}}$ the $K_{U_x}$-invariant sub-bundle
of $\widetilde{E}_{U_x}$ on $\widetilde{U}_x$,
then $(G_{U_x}^{E}/K_{U_x},\widetilde{E}_{U_x}^{\mathrm{pr}})$ defines
an orbifold sub-bundle $E^{\mathrm{pr}}$ of $E$ on $X$.
We call $E^{\mathrm{pr}}$ the proper part of $E$.
 We say $E$ is proper
if $G_{U_x}^{E}= G_x$ for any $x\in X$.
%For any orbifold vector bundle $E$, its proper part $E^{\mathrm{pr}}$
%is a proper orbifold vector bundle.

Now, any structure on $X$ or $E$ should be locally
$G_x$ or $G_{U_x}^{E}$ equivariant.

Assume that there exists a proper orbifold Hermitian line bundle
$L$ over $X$ endowed with
a Hermitian connection $\nabla^L$ with the property that
$\frac{\sqrt{-1}}{2\pi}R^L=\omega$ (Thus there exist $k\in \bN$
such that $L^k$ is a line bundle in the usual sense).
Let $(E,h^E)$ be a proper orbifold Hermitian vector
bundle  on $X$ with Hermitian connection $\nabla^E$
and its curvature $R^E$.

Then the construction in Section \ref{s2} works well here.
Especially, the spin$^c$ Dirac operator $D_p$ is well defined.
In our situation, let $\{S^p_1, \cdots, S^p_{d_p}\}$
$(d_p = \dim \Ker D_p)$ be any orthonormal basis of
$\Ker D_p$ with respect to the inner
 product (\ref{b3}). We still have (\ref{c1}) for $B_p(x)$. In fact,
 on the local coordinate above, $\widetilde{S}^p_i(\widetilde{z})$
on $\widetilde{U}_x$ are $G_x$ invariant, and
\begin{align}\label{e1}
P_p (z,z') &= \sum_{i=1}^{d_p}  \widetilde{S}^p_i ( \widetilde{z})
 \otimes (\widetilde{S}^p_i(\widetilde{z'}))^*.
\end{align}

We note that if $Q: \cC^\infty(X,E)\to \cC^\infty(X,F)$ is a pseudo-differential operator of order $m$
($m < -2n-k, k\in \bN$), and $E,F$ are proper orbifold vector bundles, then the operator $Q$ has a
$\cC^k$-kernel. In fact, $Q_{U_x}$ lifts to a  pseudo-differential operator
$\widetilde{Q}_{U_x}$ on $\widetilde{U}_x$ and for
$\widetilde{Q}_{U_x}(\widetilde{z},\widetilde{z'})$  the
 $\cC^k$-kernel on $\widetilde{U}_x \times \widetilde{U}_x$ with respect to
$dv_{\widetilde{U}_x}$,
the kernel of the operator $Q_{U_x}:  \cC^{\infty}({U_x},
E_{|{U_x}})\rightarrow \cC^{\infty}({U_x}, F_{|U_x})$
is (cf. also \cite[(2.2)]{M})
\begin{eqnarray}\label{e3}
Q_{U_x}(z,z') =   \sum_{g\in G_x}(g, 1)
\widetilde{Q}_{U_x}( g^{-1} \widetilde{z},\widetilde{z'}), \quad
(z,z')\in {U_x} \times {U_x}.
\end{eqnarray}
Indeed, for $s\in  \cC^\infty(U_x,E)$ with compact support,
then $s$ is a $G_x$-invariant section of $\widetilde{E}_{U_x}$
on $\widetilde{U}_x$, by definition,
 \begin{align}\label{0e3}
 (Qs)(z)= & \int_{\wi{U}_x} \wi{Q}_{U_x} (\wi{z},\wi{z}')
s(\wi{z}') dv_{\wi{U}_x}(\wi{z}') \\
 = & \frac{1}{|G_x|} \sum_{g\in
G_x} \int_{\wi{U}_x} (g,1)\wi{Q}_{U_x} (g^{-1}\wi{z},\wi{z}')
s(\wi{z}') dv_{\wi{U}_x}(\wi{z}')\nonumber \\
 = & \int_{U_x} \sum_{g\in
G_x}(g,1)\wi{Q}_{U_x} (g^{-1}\wi{z},z') s(z') dv_{U_x}(z').\nonumber 
\end{align}

\comment{ The following analogue of Theorem \ref{t0.1} is true,
\begin{thm}\label{t4.3} There exist smooth coefficients
$b_r(x)\in \End (\Lambda (T^{*(0,1)}X)\otimes E)_x$
with $b_0=(\det {\bf J})^{1/2} I_{\bC\otimes E}$, and $b_r(x)$ are
polynomials in $R^{TX}$, $R^{\det}$,
$R^E$ (and $R^L$)  and their derivatives with order
$\leq 2r-1$  (resp. $2r$) at $x$,
such that for any $k,l\in \bN$, there exist
$C_{k,l}>0$, $N\in \bN$ such that for any $x\in X$, $p\in \bN$,
\begin{multline}\label{0e4}
\Big |\frac{1}{p^n}B_p(x)
- \sum_{r=0}^{k} b_r(x) p^{-r} \Big |_{\cC^l} \\
\leq C_{k,l}
\Big (p^{-k-1} + p^{l/2}(1+\sqrt{p}d(x,\Sigma X) )^N
e^{-C \sqrt{p} d(x,\Sigma X)}\Big ).
\end{multline}
Moreover if the orbifold $(X,\omega)$ is K\"ahler,  ${\bf J}=J$
and the proper orbifold vector bundles $E,L$ are holomorphic on $X$,
then  $b_r(x)\in \End (E)_x$ and  $b_r(x)$ are polynomials in $R^{TX}$,
$R^E$  and their derivatives with order
$\leq 2r-1$  at $x$.
\end{thm}
}

\begin{proof}[Proof of Theorem \ref{t0.4}]
At first, we have  the analogue of  Propositions \ref{0t3.0},
\be\label{0e5}
|P_p(x,x') -  F(D_p)(x,x')|_{\cC^m(X)}
\leq C_{l,m,\var} p^{-l}.
\ee
To prove  (\ref{0e5}), we work on
$\widetilde{U}_{x_i}$, and the Sobolev norm in (\ref{c11})
is summed on $\widetilde{U}_{x_i}$.

Note that on orbifold, the property of the finite propagation speed
of solutions of  hyperbolic equations still holds if we check the proof therein
\cite[\S 7.8]{CP}, \cite[\S 4.4]{T1} as pointed out in \cite{M}.
Thus for $x,x'\in X$, if $d(x, x')\geq \var$, then
$F( D_p)(x,x')= 0$, and  given $x\in X$, $F(D_p)(x,\cdot)$ only depends
on the restriction of $D_p$ to
$B^X(x,\var)$. Thus the problem on
the asymptotic expansion of $P_p (x, \cdot)$ is local.

Now, we replace $X$ by $\bR^{2n}/G_{x_0}$, and let $\wi{L},
\wi{E}$ be the $G_{x_0}$-equivariant  vector  bundles  on
$\widetilde{U}_{x_0}$ corresponding to $L,E$ on
$\widetilde{U}_{x_0}/G_{x_0}$. In particular, $G_{x_0}$ acts
linearly and effectively  on $\bR^{2n}$. We will add a superscript
$\, \wi{}\,$ to indicate the corresponding
 objects on  $\bR^{2n}$.

Now for $Z,Z' \in \bR^{2n}/G_{x_0}$, $|Z|, |Z'|\leq \var/2$ and $
\wi{Z},\wi{Z}' \in \bR^{2n}$ represent $Z,Z'$,
then by (\ref{1c3}), (\ref{1c19}), and (\ref{e3}), for any $l,m\in \bN$,
there exists $C_{l,m,\var}>0$ such that for $p\geq 1$,
\begin{align}\label{e6}
&F(D_p)(Z,Z') = \sum_{g\in G_x}(g, 1) F(\wi{D}_p)( g^{-1} \wi{Z},\wi{Z}'),\\
&|F(\wi{D}_p)(\wi{Z},\wi{Z}')-\wi{P}^0_p(  \wi{Z},\wi{Z}')|_{\cC^m}
\leq C_{l,m,\var} p^{-l}.\nonumber
\end{align}
Moreover, for $t=\frac{1}{\sqrt{p}}$,
\begin{align}\label{0e6}
\wi{P}^0_p(\wi{Z},\wi{Z}')=\frac{1}{t^{2n}}
\wi{P}_{0,t}(\wi{Z}/t,\wi{Z}'/t)\kappa ^{-1} (Z') .
\end{align}
We will denote $P^{(r)}$ in (\ref{ue52}) by $P^{(r)}_{x_0}$
to indicate the base point $x_0$. For $g\in G_{x_0}$, we denote by
$\wi{Z}= \wi{Z}_{1,g} + \wi{Z}_{2,g}$ with $\wi{Z}_{1,g}\in T\wi{U}^g_{x_0}$,
$ \wi{Z}_{2,g}\in N_{g,x_0}$ (here  $N_{g,x_0}$ is the normal bundle
to $\wi{U}^g_{x_0}$ in $\wi{U}_{x_0}$).
 By (\ref{1ue32}), (\ref{0ue45}), as in (\ref{ue66}),
for $|\wi{Z}|\leq \var/2$, $\alpha, \alpha'$ with
$|\alpha|\leq m, |\alpha'|\leq m'$,
\begin{multline}\label{0e7}
\Big |\frac{\partial^{|\alpha|}}
{\partial \wi{Z}_{1,g}^{\alpha}}\frac{\partial^{|\alpha'|}}
{\partial \wi{Z}_{2,g}^{\alpha'}} \Big (\frac{1}{p^n} \wi{P}_p^0 (g^{-1}\wi{Z},\wi{Z})
- \sum_{r=0}^k t^r P^{(r)}_{\wi{Z}_{1,g}}
(\sqrt{p}g^{-1} \wi{Z}_{2,g}, \sqrt{p} \wi{Z}_{2,g})
\kappa ^{-1}_{\wi{Z}_{1,g}} (\wi{Z}_{2,g}) \Big )\Big | \\
\leq C t^{k-m'} (1+\sqrt{p}|\wi{Z}_{2,g}|)^N
\exp (- \sqrt{C''\mu_0} \sqrt{p}|\wi{Z}_{2,g}|).
\end{multline}
Especially, for $Z \in \bR^{2n}/G_{x_0}$, $|\wi{Z}|\leq \var/2$, as in (\ref{ue67}),
\begin{align}\label{0e8}
\sup_{|\alpha|\leq m}\Big |\frac{\partial^{|\alpha|}}
{\partial \wi{Z}^{\alpha}}\Big  (\frac{1}{p^n} \wi{P}_p^0(\wi{Z},\wi{Z})
- \sum_{r=0}^k  p^{-r} P^{(2r)}_{\wi{Z}} (0,0)\Big )\Big |\leq C p^{-k-1}.
\end{align}
Thus from (\ref{0e5})-(\ref{0e8})
\footnote{In the same way,
by Theorem \ref{tue17}, \eqref{e6}, \eqref{0e6},
we get the full off-diagonal expansion of the Bergman kernel
on the orbifolds as in Theorem 3.18$^\prime$.},
we get for $|\wi{Z}|\leq \var/2$,
\begin{multline}\label{0e9}
\sup_{|\alpha|\leq m'}\left |\frac{\partial^{|\alpha|}}
{\partial \wi{Z}^{\alpha}}\Big (\frac{1}{p^n} P_p(\wi{Z},\wi{Z})
- \sum_{r=0}^k b_r(\wi{Z})p^{-r}\right.\\
 \left. - \sum_{r=0}^{2k}  p^{-\frac{r}{2}}
\sum_{1\neq g\in G_{x_0}} (g,1)  P^{(r)}_{\wi{Z}_{1,g}}
(\sqrt{p}g^{-1} \wi{Z}_{2,g}, \sqrt{p} \wi{Z}_{2,g})
\kappa ^{-1}_{\wi{Z}_{1,g}} (\wi{Z}_{2,g})\Big)\right |\\
\leq C \left(p^{-k-1} + p^{-k+\frac{m'-1}{2}} \left(1+\sqrt{p}d (Z,X')\right)^N
\exp \left(- \sqrt{C''\mu_0} \sqrt{p}d (Z,X')\right)\right).
\end{multline}

By  (\ref{1ue59}), we get for $\alpha, \alpha'$ with
$|\alpha|\leq m, |\alpha'|\leq m'$,
\begin{align}\label{0e10}
\Big |\frac{\partial^{|\alpha|}}
{\partial \wi{Z}_{1,g}^{\alpha}}\frac{\partial^{|\alpha'|}}
{\partial \wi{Z}_{2,g}^{\alpha'}} \sum_{r=0}^{2k} t^r P^{(r)}_{\wi{Z}_{1,g}}
(g^{-1} \wi{Z}_{2,g}/t,  \wi{Z}_{2,g}/t)\Big |
\leq C t ^{-m'} (1+|\frac{\wi{Z}_{2,g}}{t}|)^N
\exp (-\frac{C'}{t} |\wi{Z}_{2,g}|).
\end{align}
For any compact set $K\subset X\setminus X'$, we get the uniform estimate
 (\ref{0.7}) from (\ref{0e8}) as in Section \ref{s3.6} as $G_x=\{1\}$.
 From (\ref{0e9}), (\ref{0e10}), we get (\ref{0.7})
near the singular set $X'$.

By the argument in Section \ref{s4.1}, we have established the last
part of Theorem \ref{t0.4}.
\end{proof}

Note that if $x_0\in X'$, then $|G_{x_0}|> 1$. Now, if in addition,
$L$ and $E$ are usual vector bundles, i.e. $G_{x_0}$ acts on both
$L_{x_0}$ and $E_{x_0}$ as identity, then by (\ref{0e9}),
\begin{align}\label{e10}
\Big |\frac{1}{p^n}P_p (x_0,x_0)
- |G_{x_0}|b_0(x_0)\Big |
\leq C p^{-1/2}.
\end{align}
Thus we can never have an uniform asymptotic expansion on $X$
if $X'$ is not empty.

\begin{rem} On $\wi{U}^g_{x_0}$, $g$ acts on $L$ by multiplication by
 $e ^{i\theta}$, the action of $g$ on $E_{\wi{U}^g_{x_0}}$ and
on $\Lambda (T^{*(0,1)}X)$ is parallel with respect to the connections
$\nabla ^E$ and $\nabla ^{\rm Cliff}$. We denote by $g|_{\Lambda \otimes E}$,
$g|_E$ the action of $g$ on $ \Lambda (T^{*(0,1)}X)\otimes E$, $E$
on  $\wi{U}^g_{x_0}$. We define on $\wi{U}^g_{x_0}$
\begin{align}\label{e11}
\psi_{r,q}(\wi{Z}_{1,g})
= \sum_{|\alpha|=q} \frac{1}{\alpha!}\left[\int_{ N_{g,x_0} }
 g|_{\Lambda\otimes E}
P^{(r)}_{\wi{Z}_{1,g}} (g ^{-1} \wi{Z}_{2,g}, \wi{Z}_{2,g})
\wi{Z}_{2,g}^\alpha dv_N(\wi{Z}_{2,g})\right]
\Big(\frac{\partial}{\partial \wi{Z}_{2,g}}\Big)^\alpha
(\kappa ^{-1}_{\wi{Z}_{1,g}}\cdot) .
\end{align}
Then $e ^{i\theta p}\psi_{r,q}(\wi{Z}_{1,g})$
are a  family of differential operators on $\wi{U}^g_{x_0}$
along the normal direction $N_{g,{x_0}}$
with coefficients in $\End (\Lambda (T^{*(0,1)}X)\otimes E)$, and
they are well defined on $\wi{U}^g_{x_0}/Z_{G_{x_0}}(g)$ and on $\Sigma X$. By (\ref{1ue59}), (\ref{0e9}), we know that in
the sense
of distributions,
\begin{align}\label{e12}
\frac{1}{p^n} B_p(x)= \sum_{r=0}^k p^{-r/2} \sum_{X_j\subset X\cup \Sigma X}\frac{1}{m_j} p^{-n+\dim X_j}e ^{i\theta_j
p}\delta_{X_j} \sum_{q\geq 0} p^{-\frac{q}{2}}\psi_{r,q}
+ \cO(p^{-k}).
\end{align}
Here $X_j$ runs over all the connected component of $X\cup \Sigma X$ and $g$ acts on $L|_{X_j}$ as multiplication by $e
^{i\theta_j}$, and $m_j$ is the multiplicity of $X_j$ defined in \cite{K} (cf. also \cite{M}).

Especially, if $\Sigma X=\{y_j\}$ is finite points, then $m_j= |G_{y_j}|$ and
$g|_{\Lambda (T^{*(0,1)}X)\otimes E}\circ I_{\bC\otimes E}= g|_E\circ I_{\bC\otimes E}$. Moreover, as $g$ commutes with
$\mathcal{J}_{x_0}$,
 from (\ref{ue62}),  for $Z=z+\overline{z}$,
\begin{multline}\label{e15}
\int_{\bR^{2n}} P (g^{-1}Z,Z) dZ \\
=
\frac{\det_{\bC}{\mathcal{J}_{x_0}}}{(2\pi)^n} \int_{\bR^{2n}}
\exp\Big (- \frac{1}{4} ||\mathcal{J}_{x_0}|^{\frac{1}{2}}(g^{-1}-1)Z|^2
+\frac{1}{2}    \left \langle \mathcal{J}_{x_0} g^{-1}Z,Z \right \rangle\Big )dZ\\
=\frac{\det_{\bC}{\mathcal{J}_{x_0}}}{(2\pi)^n} \int_{\bR^{2n}}
\exp\Big (- \frac{1}{2}\left \langle |\mathcal{J}_{x_0}|Z,Z \right \rangle
+ \frac{1}{2}\left \langle (|\mathcal{J}_{x_0}|+\mathcal{J}_{x_0})g^{-1}Z,Z \right \rangle \Big )dZ\\
= \frac{\det_{\bC}{\mathcal{J}_{x_0}}}{(2\pi)^n} \int_{\bR^{2n}}
\exp\Big (- \left \langle \mathcal{J}_{x_0} z,\overline{z}\right \rangle
+ \left \langle \mathcal{J}_{x_0}g^{-1}z,\overline{z}\right \rangle\Big )dZ\\
= \frac{1}{\det_\bC (1-g^{-1} _{T^{(1,0)}X})}.
\end{multline}
Thus  from (\ref{1ue59}), (\ref{e12}), (\ref{e15}),
\begin{align}\label{e16}
 B_p(x)= \sum_{r=0}^n b_r(x) p^{n-r} + \sum_{y_j} \frac{ e ^{i\theta_j p}g|_E\circ I_{\bC\otimes E}}{|G_{y_j}|\det_\bC
(1-g^{-1} _{T^{(1,0)}X})}\delta_{y_j}
+  \cO(\frac{1}{p}).
\end{align}
\end{rem}

\begin{rem}
Now assume that $(X,\omega)$ is a K\"ahler orbifold and ${\bf J}=J$,
moreover $L$ is an usual line bundle on $X$.
Then we can
 embed $X$ into $\field{P}(H^0(X,L^p)^*)$ by using the orbifold
 Kodaira embedding $\phi_p$
for $p$ large enough (cf. \cite[\S 7]{Baily57}).
Let $\mO(1)$ be the canonical line bundle on $\field{P}(H^0(X,L^p)^*)$
with canonical metric $h^{\mO(1)}$. Then $L^p= \phi_p^* \mO(1)$
and $h^{L^p}= B_p(x) \phi_p^*h^{\mO(1)}$. We can also interpret as following:
Let $\{S_j\}_{j=1}^{d_p}$ be
an orthonormal basis of $H^0(X,L^p)$ with respect to (\ref{b3}),
then it induces an identification
 $H^0(X,L^p)^*\cong \bC^{d_p}$;
also, choose a local $G_x$-invariant holomorphic frame $S_L$
(which is possible as $G_x$ acts on $L_x$ as identity)
and write $S_j=f_jS_L^p$.
Then $\phi_p: X\hookrightarrow \field{P}(H^0(X,L^p)^*)$
is defined by $\phi_p(x) = [f_1(x),\cdots, f_{d_p}(x)]$.
 Let $\omega_{FS}$ be the Fubini-Study metric on $\field{P}(H^0(X,L^p)^*)$.
Then
\begin{align}\label{e17}
\frac{1}{p}\phi_p^* \omega_{FS}
= \frac{\sqrt{-1}}{2\pi p}
\partial \overline{\partial} \log \Big(\sum_{j=1}^{d_p}|f_j|^2\Big)
=  \omega
+ \frac{\sqrt{-1}}{2\pi p} \partial \overline{\partial} \log B_p(x).
\end{align}
Note that $g\in G_{x_0}$ acts as identity on $L_{x_0}$,
for $\wi{Z}=z+\overline{z}$, by (\ref{ue62}),
\begin{multline}\label{e18}
(g,1)P_{\wi{Z}_{1,g}}(\sqrt{p}g^{-1}\wi{Z}_{2,g},\sqrt{p}\wi{Z}_{2,g})
\\
= \exp\Big (- \frac{\pi }{2}p \left|(g^{-1}-1)\wi{Z}_{2,g}\right|^2
+\pi p   \left \langle  g^{-1}(z_{2,g}-\overline{z}_{2,g}),\wi{Z}_{2,g} \right \rangle\Big ).
\end{multline}
Set $\wi{b}_0(\wi{Z}) = 1+ \sum_{1\neq g\in G_{x_0}}(g,1)
P_{\wi{Z}_{1,g}}(\sqrt{p}g^{-1} \wi{Z}_{2,g},\sqrt{p} \wi{Z}_{2,g})
\kappa_{\wi{Z}_{1,g}}(\wi{Z}_{2,g}) $. Then $\wi{b}_0(\wi{Z})$
has a positive real part on $T_{x_0}X$.
By (\ref{0e9}), for $m\in \bN$, taking $k\gg m$, then for $p$ large enough, for
$|Z|\leq \var/2$, under the norms $\cC^m$,
\begin{multline}\label{e19}
\log (\frac{1}{p^n}B_p(Z))= \log(\wi{b}_0(\wi{Z})) +
\log \left(1+\sum_{r=1}^k \wi{b}_0(\wi{Z})^{-1} b_r(\wi{Z})p^{-r}\right. \\
 \left. - \sum_{r=1}^{2k}  p^{-\frac{r}{2}} \wi{b}_0(\wi{Z})^{-1}
\sum_{1\neq g\in G_{x_0}} (g,1)  P^{(r)}_{\wi{Z}_{1,g}}
(\sqrt{p}g^{-1} \wi{Z}_{2,g}, \sqrt{p} \wi{Z}_{2,g})
\kappa_{\wi{Z}_{1,g}}(\wi{Z}_{2,g})\right)
+ \cO(p^{-k+\frac{m}{2}}).
\end{multline}
Thus from  (\ref{1ue59}), (\ref{e17}), (\ref{e19}),
for any $l\in \bN$, there exists $C_l>0$ such that
\begin{align}\label{e20}
\left|\frac{1}{p}\phi_p^* \omega_{FS}(x) -\omega(x)\right|_{\cC^l}\leq C_l
\left(\frac{1}{p}+ p^{\frac{l}{2}} e ^{-c\sqrt{p} d(x,X')}\right).
\end{align}
\end{rem}

%%%%%%%%%%%%%%%%%%%%%%%%%%%%%%%%%%%%%%%%%%%%%%%%%%%%%%%%%%%%%%%%%%%%%


\begin{thebibliography}{99}

\bibitem{Baily57} W. L. Baily, On the imbedding of $V$-manifolds
in projective space.  {\em Amer. J. Math.}  79  (1957), 403--430.

\bibitem{BFG} M. Beals, C. Fefferman, R. Grossman,
Strictly pseudoconvex domains in $ C\sp{n}$.
{\em Bull. Amer. Math. Soc. (N.S.)} 8 (1983), no. 2, 125--322.

\bibitem{BeGeV}
N. Berline, E. Getzler and M. Vergne,
{\it Heat kernels and Dirac operators}, Springer--Verlag, 1992.

\bibitem{B} J.-M. Bismut, Demailly's asymptotic Morse inequalities: a heat equation proof. {\em J. Funct. Anal.} 72 (1987),
no. 2,
263--278.

 \bibitem{B90} J.-M. Bismut, Koszul complexes, harmonic oscillators and the Todd class,
{\em J.A.M.S.} 3 (1990), 159-256.

\bibitem{B95}  J.-M. Bismut, Equivariant immersions and Quillen metrics,
{\em J. Diff. Geom.} 41 (1995). 53-159.


\bibitem{BL} J.-M. Bismut  and  G. Lebeau, Complex immersions and Quillen
metrics,  {\em Publ. Math. IHES.}, Vol. 74, 1991, 1-297.



\bibitem{BiV}
J.--M. Bismut and E. Vasserot,
The asymptotics of the Ray--Singer analytic torsion associated with high powers of a positive line bundle,
{\em Commun. Math. Phys.}, 125 (1989), 355--367.


\bibitem{BU}
D. Borthwick and A. Uribe,
Almost complex structures and geometric quantization,
{\em Math. Res. Lett.}, 3 (1996), 845--861. Erratum:  5 (1998), 211-212.

\bibitem{BU1}
D.~Borthwick and A.~Uribe,
 {Nearly K{\"a}hlerian Embeddings of Symplectic Manifolds},
 {\em  Asian J. Math.} 4 (2000), no.~3, 599--620.

\bibitem{Bou} T. Bouche,
Convergence de la m\'etrique de Fubini-Study d'un fibr\'e
lin\'eaire positif. {\em Ann. Inst. Fourier.} 40  (1990), 117--130.

\bibitem{BoG}
L.~{B}outet~de Monvel and V.~Guillemin, \emph{{The spectral theory of Toeplitz
  operators}}, Annals of Math. Studies, no.~99, Princeton Univ. Press,
  Princeton, NJ, 1981.

\bibitem{BS}  L. Boutet de Monvel and J. Sj\"ostrand,
Sur la singularit\'e des noyaux de Bergman et de Szeg\"o,
Ast\'erisque 34-35 (1976), 123--164.

\bibitem{CP} J. Charazain, A. Piriou,  {\em Introduction \`a la th\'eorie des
 \'equations aux d\'eriv\'ees partielles,} Paris: Gauthier-villars 1981.

\bibitem{CGT} J. Cheeger, M. Gromov, M. Taylor, Finite propagation
speed, kernel estimates for functions of the Laplace operator, and
the geometry of complete Riemannian manifolds, {\em  J.
Differential Geom.}  17 (1982), 15-53.

\bibitem{Catlin} D.  Catlin,
The Bergman kernel and a theorem of Tian.
Analysis and geometry in several complex variables (Katata, 1997), 1--23,
Trends Math.,
Birkh\"auser Boston, Boston, MA, 1999.

\bibitem{Charles} L. Charles,  Berezin-Toeplitz operators,
a semi-classical approach.   {\em  Comm. Math. Phys.}  239  (2003), 1--28.

\bibitem{Chernoff} P.R. Chernoff,
Essential self-adjointness of powers of generators of hyperbolic equations.
J. Functional Analysis 12 (1973), 401--414.


\bibitem{Christ91} M. Christ,
On the $\overline{\partial}$ equation in weighted
$L^2$-norms in $\bC ^1$, {\em J. Geom. Anal.} 3 (1991), 193-230.

\bibitem{DKM}  X. Dai, K. Liu, X. Ma,
On the asymptotic expansion of Bergman kernel.
{\em C.R.A.S. Paris}, 339 (2004), 193-198.


\bibitem{Demailly} J.P. Demailly, Holomorphic Morse inequalities.
Several complex variables and complex geometry, Part 2 (Santa Cruz, CA, 1989),
 93--114, Proc. Sympos. Pure Math., 52,
Part 2, Amer. Math. Soc., Providence, RI, 1991.

\bibitem{D} S. K. Donaldson,  Scalar curvature and projective
embeddings. I. {\em  J. Differential Geom.} 59 (2001), no. 3, 479--522.

\bibitem{K} T. Kawasaki, The Riemann-Roch theorem for  V-manifolds.
{\em  Osaka J. Math} 16 (1979) 151-159.

\bibitem{LM}
H. B. Lawson and M.-L. Michelson, {\it Spin geometry}, Princeton Univ. Press,
Princeton, New Jersey, 1989.

\bibitem{Li} N. Lindholm,  Sampling in weighted $L\sp p$ spaces of entire
functions in $\bC\sp n$ and estimates of the Bergman kernel.
{\em  J. Funct. Anal.} 182 (2001), 390--426

\bibitem{Lu} Z. Lu,  On the lower order terms of the asymptotic
 expansion of Tian-Yau-Zelditch. {\em  Amer. J. Math.}  122 (2000), no. 2, 235--273.

\bibitem{M} X. Ma, Orbifolds and analytic torsions.
{\em  Trans. Amer. Math. Soc.} 357 (2005), 2205-2233.

\bibitem{MM} X. Ma, G. Marinescu,
 The ${\rm spin}\sp c$ Dirac operator on high tensor powers of a
line bundle. {\em  Math. Z.}  240 (2002), no. 3, 651--664.

\bibitem{Ru} W. Ruan, Canonical coordinates and Bergman metrics.
{\em  Comm. Anal. Geom.}  6  (1998), 589--631.

\bibitem{SZ02}
B.~Shiffman and S.~Zelditch, Asymptotics of almost holomorphic sections
  of ample line bundles on symplectic manifolds, {\em J. Reine Angew. Math.}
  544 (2002), 181--222.

\bibitem{T1} M. Taylor, {\em Pseudodifferential operators},
Princeton Univ Press, Princeton 1981.


\bibitem{Tian}   G. Tian, On a set of polarized K\"ahler metrics on
algebraic manifolds, J. Differential Geom. 32 (1990), 99--130.

\bibitem{Wang1} X. Wang,  Thesis, 2002.

\bibitem{Zelditch}  S. Zelditch,
Szeg\"o kernels and a theorem of Tian.
{\em  Internat. Math. Res. Notices} 1998, no. 6, 317--331.

\end{thebibliography}
\end{document}